\documentclass{fic-l}

  \newtheorem{theorem}{Theorem}[section]

\theoremstyle{definition}

\newtheorem{example}[theorem]{Example}

\theoremstyle{remark}

\numberwithin{equation}{section}

\usepackage{amscd}

\usepackage{graphicx}
\usepackage{wrapfig}
 \newtheorem{proposition}[theorem]{Proposition}

\title{Symplectic Gluing and Family Gromov-Witten Invariants}

\author{ Junho Lee and Thomas H. Parker}
\thanks{The second author was partially supported
by the N.S.F.}
\address{Michigan State University\\ East Lansing,
MI   48824}

\date{\today}
\addtocounter{section}{0}

\begin{document}

\maketitle
 
 \begin{abstract}
This article describes the use of symplectic
cut-and-paste methods to compute Gromov-Witten invariants.  Our focus 
is on recent advances extending  these methods  to K\"{a}hler 
surfaces with geometric genus $p_g>0$, for which the usual GW
invariants vanish for most homology classes. This involves extending the Splitting Formula and
the Symplectic Sum Formula to the family GW invariants introduced by
the first author.   We present applications to the invariants
of elliptic surfaces and to the Yau-Zaslow Conjecture.  In both cases 
the results agree with the conjectures of algebraic geometers and 
yield a proof, to appear in \cite{ll}, of previously unproved cases of 
the Yau-Zaslow Conjecture.
\end{abstract}

\vskip.4in

Gromov-Witten invariants are counts of holomorphic curves in a
symplectic manifold $X$.  To define them one  chooses an
almost  complex structure $J$ compatible with the symplectic
structure  and considers the set of maps $f:\Sigma\to X$ from
Riemann surfaces
$\Sigma$ which satisfy the (nonlinear elliptic)  $J$-holomorphic map equation
\begin{equation}
\overline{\partial}_Jf=0.
\label{1}
\end{equation}
   After compactifying the moduli space of such
maps, one  imposes constraints, counting, for example, those maps
whose images  pass through  specified points.  With the
right number of constraints and a generic perturbation of the equation, the number of such
maps is finite.  That number is a GW invariant of the symplectic
manifold $X$.

The first part of this article focuses on cut-and-paste methods for computing
GW invariants of symplectic  four-manifolds.    Two useful techniques 
are described in Sections 2--4.  The first is
the method  of ``splitting the domain'', in
which  one considers
maps whose domains are pinched Riemann surfaces.  This produces
recursion relations, called TRR formulas,  relating
one GW invariant to invariants whose images have
smaller area or lower  genus.

   One can also consider the behavior of GW
invariants as one  ``splits the target'' using the symplectic cut
operation and its inverse, the symplectic sum.
The resulting  ``Symplectic Sum Formula'' of E. Ionel and the second author
  expresses the GW invariants of the symplectic sum of two
manifolds  $X$ and $Y$ in terms of the  relative GW invariants of $X$
and $Y$ as described in Section 4.

These splitting formulas for the domain and the target can sometimes be
combined  to completely determine a set of GW invariants.  One such 
computation is given  in Section 5.   There we use a clever
geometric  argument of E. Ionel to derive an explicit formula for certain
GW invariants of the rational elliptic surface $E(1)$.  That formula,
explained in Section 5, is
\begin{equation}
\label{introFormula1}
   F(t)\,=\,\prod_{d\geq 1}
   \left(\frac{1}{1-t^{d}}\right)^{12}
\end{equation}

There is strong motivation for generalizing  this
$E(1)$ computation.  We would, of course, like to make such cut-and-paste
arguments into  a systematic way of computing invariants of whole 
classes of 4-manifolds.
More alluring yet is the  connection with enumerative algebraic
geometry.   There is a series of  conjectures, originally 
coming from
string theory, that claim that the generating functions of counts of
holomorphic curves in algebraic surfaces are given by certain
specific formulas (see \cite{go}).  The simplest such formula is exactly
(\ref{introFormula1}).  The next simplest is the Yau-Zaslow 
conjecture for rational curves in K3 surfaces (\cite{yz}).

In seeking to generalize the $E(1)$ computation, one 
encounters a serious problem.  For $E(1)$, the GW invariants are the 
same as the curve counts of enumerative geometry.  But for most other 
complex surfaces, beginning with $K3$,  the conjectural 
formulas count curves in classes for which 
the Gromov-Witten invariants  vanish!

This discrepancy occurs because GW invariants are defined using
generic almost  complex structures
$J$, while  K\"ahler structures  are very special. Entire
families of curves can disappear when the K\"ahler $J$ is perturbed
to a generic $J$.
For example,  algebraic K3 surfaces  contain  many holomorphic curves, as
predicted by the Yau-Zaslow conjecture, but a generic $K3$ surface  has no
holomorphic curves at all!  In general,
given a space (generalized Severi variety)  of curves in a 
K\"ahler manifold,  one
can perturb to a generic almost  complex structure $J$ and ask which
curves in the Severi variety perturb to become
$J$-holomorphic curves. This problem fits into standard deformation
theory, as explained in Section 6, and from that perspective one sees 
that there is an
obstruction   bundle   over the Severi variety  whose rank is the 
geometric genus $p_g=\mbox{dim }H^{0,2}(X)$ of the surface $X$.  Thus
when $p_g>0$ the space of curves and the GW moduli space have 
different dimensions, and consequently the GW invariants are 
unrelated to the enumerative counts of algebraic 
geometry.

\medskip

Sections 7--10 describe  the first author's geometric analysis  method for overcoming this 
problem.  The key observation is that on each K\"{a}hler manifold there is  a very natural family of almost complex structures $J_\alpha$ parameterized by  $H^{0,2}(X)$.  One can then consider pairs $(f,\alpha)$ where $f$ is a $J_\alpha$-holomorphic map for some $\alpha\in H^{0,2}(X)$.  Such pairs form a moduli space whose dimension is $p_g$ greater than the standard GW moduli space.  The family GW invariants are defined from that moduli space using the constructions of Li-Tian \cite{LT}. (Bryan-Leung and Behrend-Fantechi, have alternative approaches which use algebraic geometry).

The family moduli space has a remarkable property: whenever $J$ is K\"{a}hler, all $J_\alpha$-holomorphic maps in fact have $\alpha=0$.  Thus the family invariants, which are generally counts of $J_\alpha$-holomorphic maps, reduce to counts of true holomorphic maps in the K\"{a}hler surface.   A less pleasing aspect is that the moduli space may not be compact , since $\alpha$ ranges over the vector space $H^{0,2}(X)$.    However, results in \cite{l1} show that it is compact in many cases of interest,  and for those cases the family invariants are well-defined symplectic invariants that are closely related to enumerative counts.

The next task is to extend the cut-and-paste formulas to the family GW  invariants.  In Sections 8 and 9 we develop TRR and symplectic sum formulas for the family GW invariants, following \cite{l2}. Those formulas immediately allow us to extend the $E(1)$ calculation of Section 5 to the elliptic surfaces $E(n)$.  The result is a simple expression for the family GW invariants which generalizes (\ref{introFormula1}).

The final section describes the very recent work of the first author and N.C. Leung on the Yau-Zaslow conjecture.  The conjecture predicts a specific formula --- a slight modification of (\ref{introFormula1}) --- for the generating function of the counts of  rational curves in K3 surfaces.  In \cite{bl1},   Bryan and Leung proved the  conjecture in the cases when the curves represent primitive homology classes.  Sections 9 and  10 describe an independent proof of that fact and its extension to classes that are twice a primitive class. The proofs repeatedly use the symplectic cut-and-paste techniques developed in earlier sections.

\vskip1cm

%%%%%%%%%%%%%%%%%%%%%%%%%%%%%%%%%%%%%%%%%%%%%%%%%%%%%%%%%%%%
%%%%%%%%%%%%%%%%%  Section 1  %%%%%%%%%%%%%%%%%%%%%%%%%%%%%%
%%%%%%%%%%%%%%%%%%%%%%%%%%%%%%%%%%%%%%%%%%%%%%%%%%%%%%%%%%%%%%

\setcounter{equation}{0}
\section{Gromov-Witten Invariants}
\label{section1}
\bigskip

   Fix a closed symplectic four-manifold $(X,\omega)$.  Building on  ideas of
Donaldson and Gromov, one can define
symplectic invariants by introducing an almost complex structure $J$  and
counting (with orientation) the number
of $J$-holomorphic maps into $X$ satisfying certain constraints.
   This section gives an overview of the setup.

Given $(X,\omega)$, one can always choose an  almost complex
structure $J$ compatible with $\omega$, i.e. a bundle map $J:TX\to TX$
with $J^2=-\mbox{Id}$ and so that, for tangent vectors $u$ and $v$,
$g(u,v)= \omega(u,Jv)$  is a Riemannian metric.  Let ${\mathcal J}$ be
the space of such $J$.  For $J\in{\mathcal J}$, a map $f: C \to
X$ from a complex curve $C$  is  {\em $J$-holomorphic} if it satisfies
\begin{equation}
\overline{\partial}_Jf=0
\label{1.holomorphicmapeq}
\end{equation}
where $\overline{\partial}_Jf =\frac12(df +  J \circ df \circ j)$.
The $J$-holomorphic maps from a  smooth connected domain $C$ with
(distinct)  marked points $x_1, \dots , x_k$ gives a space
\begin{equation}
{\mathcal M}_{g,k}(X,A)
\label{1.modulispace}
\end{equation}
labeled by the genus $g$ of $C$, the number $k$ of marked points,
and the class $A\in H_2(X)$  that the image represents.

        The space (\ref{1.modulispace}) has a natural compactification,
which can be built in two steps. First, when $X$ is a single point
and $k\geq 3-2g$,  (\ref{1.modulispace}) is the Deligne-Mumford space
of curves.   That has a compactification
   $$
{\overline{\mathcal M}}_{g,k}
   $$
   consisting of stable connected curves with arithmetic genus $g$ and
$k$ marked points.  (A curve is {\em stable} if (i) it has only nodal
singular points, none of which are marked points, and (ii) each
irreducible component of genus $g_i$ has at least $3-2g_i$ marked or
nodal points.)  The Deligne-Mumford  compactification is an orbifold,
and is  a complex projective variety  of dimension $3g-3+k$.

   A   $J$-holomorphic map from a connected nodal curve is  {\em
stable}  if the image of  every  unstable irreducible component  of
$C$  represents a non-trivial element of   $H_2(X)$.   For technical
reasons it is helpful, following Ruan and Tian \cite{rt1}, \cite{rt2} to
replace (\ref{1.holomorphicmapeq}) by the inhomogeneous equation
\begin{equation}
\overline{\partial}_Jf=\nu
\label{1.pertholomorphicmapeq}
\end{equation}
where  $\nu$ is an appropriate perturbation term.   With that, the
space of stable maps
\begin{equation}
{\overline{\mathcal M}}_{g,k}(X,A)
\label{1.Mstable}
\end{equation}
is compact and, for generic $(J,\nu)$, is an orbifold with
\begin{equation}
\mbox{dim}\ {\overline{\mathcal M}}_{g,k}(X,A)=2\left[(n-3)(1-g)- 
K\cdot A\right] +2k
\label{1.dim}
\end{equation}
where $\dim X =2n$ and $K$ is the canonical class of $(X,J)$.  This
space has an orientation and comes with a map
\begin{equation}
st \times ev:{\overline{\mathcal M}}_{g,k}(X,A) \to 
{\overline{\mathcal M}}_{g,k} \times X^k
\label{1,evalmap}
\end{equation}
   that associates to a stable map $(f,C)$  the stabilization $st(C)$
   of its domain (obtained by collapsing all unstable irreducible components
   of $C$ to points) and the images $ev_i(f,C)=f(x_i)$ of the marked points.
    We can then pushforward the fundamental homology class of 
${\overline{\mathcal M}}_{g,k}(X,A)$
    by the map (\ref{1,evalmap}) and evaluating on cohomology classes
    $\beta\in H^*( {\overline{\mathcal M}}_{g,k}, { \Bbb Q} )$ and 
$\gamma_1,\dots, \gamma_k\in H^*(X;{ \Bbb Q})$.
     The resulting numbers are the Gromov-Witten invariants
   $$
   GW_{g,k}(X,A)(\beta; \gamma_1,\dots, \gamma_k) =
  ( \beta\cup \gamma_1 \cup \dots \cup \gamma_k)\cap \left((st \times ev)_*
   \left[{\overline{\mathcal M}}_{g,k}(X,A) \right]\right).
$$
Standard cobordism arguments  show that these are independent of
the choice of generic $(J,\nu)$, and hence are symplectic
invariants. Once we have fixed a space $X$, we will often write GW invariants
   $$ GW_{A,g}(X)(\beta; \gamma_1,\dots, \gamma_k)
      \ \ \ \mbox{as\ simply }\ \ \
      GW_{A,g}(\beta; \gamma_1,\dots, \gamma_k). $$

   The Gromov-Witten invariants are counts of the number of perturbed
   holomorphic maps that satisfy the constraints imposed by the $\beta$
   and $\gamma$ classes.  The $\gamma$ constraints have a simple geometric
   interpretation:  if one chooses pseudomanifolds $V_i\subset X$ representing
   the Poincar\'{e} dual of $\gamma_i$, then, assuming that $(J,\nu)$ is generic
   and that the $V_i$ are in general position,
   $ GW_{A,g}(\gamma_1,\dots, \gamma_k) $ is an oriented count of
   the number of $(J,\nu)$-holomorphic maps $f$ from a genus $g$ curve
   with $k$ marked points $x_1,\dots x_k$ that represent the class
   $A$ and
   satisfy $f(x_k)\in V_i$ for each $i$.

This geometric interpretation
implies several properties of the GW invariants.  For example, the
invariant vanishes unless the total degree of the constraints (the
cohomology classes $\beta$ and $\gamma_i$) is equal to the dimension
(\ref{1.dim}) of the space of stable maps.   An especially useful
property  is the ``Divisor Axiom'':  for any degree 2 constraint
$\gamma_1$ with Poincar\'{e} dual $V_1$ we have
\begin{equation}
GW_{A,g}(\gamma_1,\dots, \gamma_k)= (A\cdot V_1)\
GW_{A,g}(\gamma_2,\dots, \gamma_k).
\label{divisoraxiom}
\end{equation}
This holds because, when $V_1$ is in general position,
the image of any stable map $f:C\to X$ representing $A$
intersects
$V_1$ transversally at $A\cdot V_1$ points, so there are  exactly
$A\cdot V_1$  possibilities for the first marked point.

   The geometric meaning of the $\beta $ constraints is more subtle.
   But certain $\beta $ constraints play an important role in applications,
   as the examples in the next two sections show.

\vskip1cm

%%%%%%%%%%%%%%%%%%%%%%%%%%%%%%%%%%%%%%%%%%%%%%%%%%%%%%%%%%%%
%%%%%%%%%%%%%%%%%  Section 2  %%%%%%%%%%%%%%%%%%%%%%%%%%%%%%
%%%%%%%%%%%%%%%%%%%%%%%%%%%%%%%%%%%%%%%%%%%%%%%%%%%%%%%%%%%%%%
\setcounter{equation}{0}
\section{Splitting the Domain}
\label{section2}
\bigskip

The Deligne-Mumford space $\overline{{\mathcal M}}_{g,k}$ has  a
complex codimension-one subset whose generic element is a stable curve with one node.
Consider the case, shown below,  in which the node separates $C$ into
two components $C_1$ and $C_2$. The genus  decomposes as
$g=g_{1}+g_{2}$ and, because the node cannot be a marked point,  the
marked and nodal points separate  into two sets $x_1, \dots ,
x_{k_1+1}$ and $y_1, \dots, y_{k_2+1}$ with $k=k_{1}+k_{2}$.  There
is then a map
\vspace{-.3cm}
    \begin{center}
\begin{minipage}{2.5in}
$\sigma : \overline{{\mathcal M}}_{g_{1},k_{1}+1}\times\overline{{\mathcal M}}
_{g_{2},k_{2}+1}\to \overline{{\mathcal M}}_{g,k}$
   \end{minipage}
 \begin{minipage}{1.5in}
\includegraphics{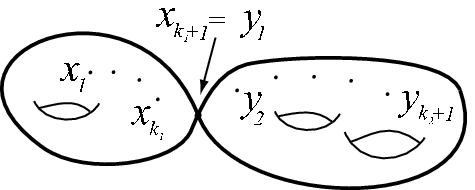}
 \hspace*{1.3cm}  {\small{\bf Figure 
1}}
     \addtocounter{figure}{1}
  \end{minipage}
  \end{center}
\vskip.4cm
\noindent defined by identifying $x_{k+1}$ with $y_{1}$. We
denote the Poincar\'{e} dual of the image of this
map $\sigma$ by  $PD(\sigma)$.   Gromov-Witten invariants with
constraint $PD(\sigma)$
\begin{equation}
\label{GWCompostion1}
GW_{A,g}(PD(\sigma );\gamma _{1},\cdots ,\gamma _{k})
\end{equation}
give counts for stable holomorphic maps from nodal domains of the type shown.

   Alternatively, a map $f:C\to X$ from such a nodal domain can also be
thought of as a pair of maps $f_1:C_1\to X$ and $f_2:C_2\to X$ with
$f_1(x_{k_1+1})=f_2(y_{1})$.  The set of such pairs $(f_1, f_2)$  can
be expressed in terms of the evaluation maps
\begin{equation}
\label{evalmap1}
ev_{k_1+1,1}: \overline{{\mathcal
M}}_{g_{1},k_{1}+1}(X,A_1)\times\overline{{\mathcal M}}
_{g_{2},k_{2}+1}(X,A_2)\to X\times X
\end{equation}
that record the images of the marked points $x_{k_1+1}$ and $y_{1}$:
the stable maps from nodal domains are those in the inverse image of
the diagonal $\Delta$ under (\ref{evalmap1}).   Thus the nodal maps
define a homology class
\begin{equation}
\label{compostion1}
\left[ \overline{{\mathcal M}}_{g_{1},k_{1}+1}(X,A_1)\right]
\otimes\left[\overline{{\mathcal M}} _{g_{2},k_{2}+1}(X,A_2)\right]
\ \cap\ \Delta^*
\quad\mbox{where} \quad
\Delta^*= ev^*_{k_1+1,1}(PD[\Delta])
\end{equation}
By the gluing theorem of Ruan-Tian \cite{rt1}, this  gives a second
description of the  Gromov-Witten invariants (\ref{GWCompostion1}).

We can go one step further.  Note that the class (\ref{compostion1})
depends only on the homology class of the diagonal in $H_{\ast
}(X\times X;{\Bbb Z})$.  We can then use the  ``splitting the
diagonal''  formula:  if $\{H_{a}\}$ and $\{H^{a}\}$ are bases of
$H^{\ast }(X;{\Bbb Z})$ dual by the intersection form, we have
\begin{equation}
\label{splitdiag}
PD[\Delta] = \sum_a H^{a} \otimes H_{a}.
\end{equation}
Inserting this in (\ref{compostion1}) and evaluating on classes
$\gamma _{1},\cdots ,\gamma_k$ then yields the first splitting
formula.

\begin{proposition}[Splitting Formula]
\label{cl}
With $\sigma$ as in Figure 1 and any
classes $\gamma _{1},\cdots ,\gamma_k\in H^*(X;{ \Bbb Q})$, we have

\begin{eqnarray*}
   GW_{A,g}(PD(\sigma );\gamma _{1},\cdots ,\gamma _{k})\  = \ \hspace*{2in}\\
   \sum_{A=A_1+A_2}\sum_{a}\ GW_{A_1,g_{1}}(\gamma _{1},\cdots ,\gamma
_{k_{1}},H_{a})\,GW_{A_2,g_{2}}(H^{a},\gamma _{k_{1}+1},\cdots
,\gamma _{k}). \\
\end{eqnarray*}
\end{proposition}

   There is a parallel story for domains with non-separating nodes.
   A genus $g-1$ curve with $k+2$ marked points becomes a genus $g$ curve with
   $k$ marked points and a single non-separating node by identifying the last
    two marked points.  That gives a map
   \vspace*{.5cm}
    \begin{center}
\begin{minipage}{1.8in}
$\theta:\overline{{\mathcal M}}_{g-1,k+2} \to \overline{{\mathcal M}}_{g,k}$
   \end{minipage}
 \begin{minipage}{1.5in}
\includegraphics{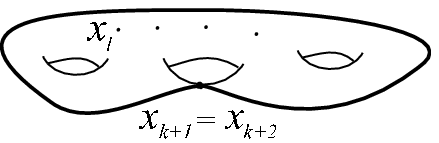}
  \end{minipage}
   \hspace*{.5cm}
 \begin{minipage}{.7in}
   \vskip.5cm
  {\small{\bf Figure 2}}
 \addtocounter{figure}{1}
    \end{minipage}
  \end{center}
\vskip.3cm
Stable maps from such genus $g$ domains are counted by GW invariants
constrained by the Poincar\'{e} dual $PD(\theta)$ of the image of
this map $\theta$.  Alternatively, they are counted by maps $f:C\to
X$ from genus $g-1$ domains with $k+2$ marked points and with
$f(x_{k+1})=f(x_{k+2})$; this  constraint is the pullback of
$PD(\Delta)$ under the evaluation map
$$
ev_{k+1,k+2}:\overline{{\mathcal M}}_{g-1,k+2}(X,A) \to X \times X
$$
   analogous to (\ref{compostion1}). Again splitting the diagonal in
homology, we get the splitting formula for non-separating nodes:
\begin{eqnarray}
\label{splittingformula2}
    GW_{A,g}(PD(\theta );\gamma _{1},\cdots ,\gamma_{k})
    & = & \sum_{a}GW_{A,g-1}(\gamma _{1},\cdots
,\gamma_{k},H_{a},H^{a}) \nonumber \\
     & = & (A\cdot A) \ GW_{A,g-1}(\gamma _{1},\cdots ,\gamma_{k})
\end{eqnarray}
    where the last equality holds on 4-manifolds.
    \medskip

    The splitting formulas   relate GW invariants in a given genus to
invariants with smaller symplectic area or lower  genus.  We will see later how  such
relations can be encoded as differential equations on generating
functions for GW invariants and used to help determine those
invariants.

\vskip 1cm

%%%%%%%%%%%%%%%%%%%%%%%%%%%%%%%%%%%%%%%%%%%%%%%%%%%%%%%%%%%%
%%%%%%%%%%%%%%%%%  Section 3  %%%%%%%%%%%%%%%%%%%%%%%%%%%%%%
%%%%%%%%%%%%%%%%%%%%%%%%%%%%%%%%%%%%%%%%%%%%%%%%%%%%%%%%%%%%%%
\setcounter{equation}{0}
\section{TRR and Descendents}
\label{section3}
\bigskip

 GW invariants
constrained by cohomology classes in ${\overline{\mathcal M}}_{g,k}$ can sometimes be rewritten and simplified using the collection of formula called the 
``Topological Recussion Relations'' (TRR).   TRR
formulas will play an important role in later sections.  The basic idea is the same as in the previous section:   relations among cohomology classes in Deligne-Mumford spaces give rise to relations among GW invariants.   Here we
will explain one particular  TRR formula, which involves  genus 1 GW invariants
constrained by ``descendent classes''.  We also describe the
geometric meaning of those constraints in simple cases.  The article
\cite{ge} of E. Getzler is a good reference for this material.

This genus 1 TRR  formula originates from  a fact about the space
${\overline{\mathcal M}}_{1,1}$ of complex structures on a torus with 
one marked point.
This is an orbifold whose boundary ${\overline{\mathcal 
M}}_{1,1}\setminus{\mathcal M}_{1,1}$
consists of a single point, which corresponds to a
curve of geometric genus $g=0$ with one non-separating node.  There
is a natural line bundle ${\mathcal L}\to {\overline{\mathcal 
M}}_{1,1}$ whose fiber at a
curve $(C,x)$ is the cotangent space
$T_{x}^{*}C$.  The chern class of $c_1({\mathcal L})$, which is
usually denoted by $\phi_1$,   is a multiple of the Poincar\'{e} dual
of
the orbifold fundamental class $\delta$ of the boundary:

\begin{equation}\label{trr1}
      \phi_{1}\ \overset{def}= c_{1}({\mathcal L})\ =\
      \frac{1}{12}\,\delta
      \mbox{\ \ \ \ in\ \ \ \ }
      H^{2}(\,{\overline{\mathcal M}}_{1,1};{\Bbb Q}\,).
\end{equation}

There is a corresponding story on the space ${\overline{\mathcal 
M}}_{1,1}(X,A)$ of
stable maps into $X$.  There is again a line bundle
$\overline{{\mathcal L}}$ whose fiber at $(C,x,f)$ is the cotangent
space
$T_{x}^{*}C$; its first chern class is denoted by $\psi_1$.  By the
definition of stable map, each irreducible component of the domain is
either a stable curve or has non-trivial image in $X$.  Under the
stabilization map
$$
st: {\overline{\mathcal M}}_{1,1}(X,A)\to {\overline{\mathcal M}}_{1,1}
$$
of (\ref{1,evalmap}),  the inverse image of $\delta$ consists of
maps from a rational nodal curve, possibly with one or more unstable
$g=0$ ``bubble'' components.  Over the open dense set of
${\overline{\mathcal M}}_{1,1}(X,A)$ consisting of maps with 
irreducible domains, the
relative cotangent bundles are related by $\overline{{\mathcal
L}}=st^*{\mathcal L}$, but that is not the case for maps with
unstable components --- there are correction terms.
The correct relation can be written symbolically as
\begin{equation}\label{trr}
      \psi_1\ =\ st^*\phi_1 \ +\
\begin{minipage}{.47in}
\includegraphics{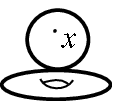}
\end{minipage}
\ =\ \frac{1}{12}\
\begin{minipage}{.45in}
\includegraphics{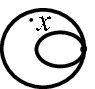}
\end{minipage}
    +\
\begin{minipage}{.47in}
\includegraphics{LP_TRR.eps}
\end{minipage}
\end{equation}
This equation requires some interpretation.  The pictures represent
the (real) codimension two  stratum of ${\overline{\mathcal 
M}}_{1,1}(X,A)$ that is  the
closure of the  set of   all stable maps whose domains have the form
shown; in particular the ``pinched torus'' in the middle picture
stands for   $st^{-1}(\delta)$.   With that understood, the symbolic
formula (\ref{trr}) says that the  GW invariant with constraint
$\psi_1$ is equal to $\frac{1}{12}$ the invariant  with constraint
$\delta$ plus the  contribution of the stratum  whose generic element consists of stable
maps from a smooth torus with one $g=0$ bubble attached.

   In practice, the GW invariants constrained by $\delta$ (the middle
picture)  can be  expressed in terms of genus 0  GW invariants  using
the Splitting Formula of GW invariants.  The invariants constained as
in the last picture can similarly be expressed in terms of invariants
corresponding to the two components.

   \bigskip

We can use $\psi_1$ to impose constraints on GW invariants.   That is
commonly done by using the evaluation map $ev: {\overline{\mathcal 
M}}_{1,1}(X,A)\to X$
(which simply evaluates each stable map at the given marked point) to
introduce ``descendent classes''
\begin{equation}\label{descendent}
\tau(\gamma)= \psi_1 \cup ev^*\gamma
\end{equation}
for each cohomology class $\gamma\in H^*(X,{\Bbb Q})$.  With more marked
points, we have a similar evaluation map $ev_i: {\overline{\mathcal 
M}}_{g,k}(X,A)\to X$
for each marked point and line bundles $\overline{{\mathcal L}}_i$
(whose fiber at $(f,C,x_1, \dots, x_k)$ is the cotangent space
$T^*_{x_i}C$) and corresponding descendent classes $\tau(\gamma_i)=
\psi_i \cup ev_i^*\gamma_i$ where $\psi_i=c_1(\overline{{\mathcal
L}}_i)$.  We then get constrained GW invariants by evaluating
products of such descendent classes on the fundamental homology class
of the space of stable maps. The simplest example occurs when
$2+\mbox{deg}\ \gamma $ is equal to the dimension of the moduli
space; then the constrained invariant with one descendent
\begin{equation}\label{descendentGW}
GW_{g, 1}(X,A)(\tau(\gamma))\ =\ \left(\psi_1 \cup
ev^*\gamma\right)\,\cap\,\left[ {\overline{\mathcal M}}_{g,1}(X,A)\right]
\end{equation}
is a number.

   This invariant (\ref{descendentGW}) has a nice  geometric
interpretation when the Poincar\'{e} dual of  $\gamma$ can be
represented by a codimension two $J$-invariant submanifold $V_\gamma\subset X$ with
trivial normal bundle (as happens in the computation of the next
section). In that case we can  fix a trivialization of the normal
bundle to $V_\gamma$ and make two observations:

\bigskip

\noindent $\bullet$\ \  The ``cutdown moduli space''
${\overline{\mathcal M}}(\gamma)\subset {\overline{\mathcal 
M}}_{g,1}(X,A)$ consisting of all maps
   $f:(C,x)\to X$ which take the marked point $x$ onto $V_\gamma$
represents $ ev^*\gamma \cap \left[ {\overline{\mathcal M}}_{g,1}(X,A)\right]$.

\medskip

\noindent $\bullet$\ \   For each map $f:(C,x)\to X$, the
normal component $df^N_x$ of the differential at the marked point $x$
is an element  of $T^*_xC$.  Thus $f\mapsto df^N_x$ is a section of
$\overline{{\mathcal L}}_1\to {\overline{\mathcal M}}_{g,1}(X,A)$. 
The zeros of this
section along ${\overline{\mathcal M}}(\gamma)$ represent the invariant
(\ref{descendentGW}) {\em provided those zeros occur only at maps in
${\overline{\mathcal M}}(\gamma)$ that have no components with images 
in  $V_\gamma$}.

\bigskip

Thus, when this last italicized caveat holds, the descendent number
(\ref{descendentGW}) counts the number of stable maps $f$ with
$f(x)\in V$ and $df_x^N=0$, that is, those maps whose image contacts
$V$ to order at least two at $f(x)$.

This story extends to higher order contact.  When $f^N$ vanishes to
order $m$ at $x=x_1$, the leading coefficient in its Taylor series is
naturally an element of $Sym^m(T_xC)$.  Intrinsically, the first $m$
terms in the Taylor series is a section of the jet bundle $J^m(x)$,
and we have the ``jet exact sequence''
$$
0\to J^{m-1}(x)\to J^m(x)\to Sym^m(T_xC)\to 0.
$$
   Globally over the set of stable maps satisfying conditions a) and b)
above, the  first $m$ terms in the Taylor series  give a section of a
relative jet bundle ${\mathcal J}^m$ whose Euler class, computed from
the jet exact sequence and induction, is
$$
c_1(\overline{{\mathcal L}_1}^m)\cup c_1(\overline{{\mathcal
L}_1}^{m-1})\cup\ \dots\ \cup c_1(\overline{{\mathcal L}_1}) \cup
ev^*_1(\gamma)
\ =\  m!\ \psi_1^m\cup ev_1^*(\gamma).
$$
Thus it is natural to consider not just $\tau$, but also the more
general descendent classes  $\tau_m$ defined by
$$
\tau_m(\gamma_i)=\psi_i^m\cup ev^*_i(\gamma_i)
$$
and the corresponding descendent invariants
$$
GW_{g, k}(X,A)\left(\tau_{m_1}(\gamma_1), \ \dots\ ,
\tau_{m_k}(\gamma_k)\right)
$$
which one sees in the literature.   In the simplest case of a single
marked point, the discussion above shows that the invariant
$$
GW_{g, 1}(X,A)\left(\tau_{m}(\gamma) \right)
$$
has a straightforward interpretation whenever conditions in a) and b)
above hold:  it   counts the number of stable maps $f$ whose image
contacts $V_\gamma$ to order $m+1$ at the point $f(x)$.  We will use
this geometric interpretation (in the case $m=1$) in Section 5 below.

\vskip.8cm
%%%%%%%%%%%%%%%%%%%%%%%%%%%%%%%%%%%%%%%%%%%%%%%%%%%%%%%%%%%%
%%%%%%%%%%%%%%%%%  Section 4  %%%%%%%%%%%%%%%%%%%%%%%%%%%%%%
%%%%%%%%%%%%%%%%%%%%%%%%%%%%%%%%%%%%%%%%%%%%%%%%%%%%%%%%%%%%%%

\setcounter{equation}{0}
\section{Splitting the Target:  The Symplectic Sum Formula}
\label{section4}
\bigskip

Another  powerful compuational tool is the   Symplectic Sum Formula.
This  can be viewed as a fact about maps into a symplectic fibration.
Consider a fibration $\lambda : {\mathcal Z} \to D$ over the disk 
whose fibers are smooth
symplectic manifolds for     $\lambda \neq 0$. 
\begin{wrapfigure}{l}{4.4cm}
\begin{minipage}{4.4cm}
\includegraphics[scale=1.5]{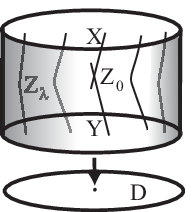} 
\caption{}
\end{minipage}
\end{wrapfigure}

 Suppose
that the central fiber $Z_{0}$ is the singular symplectic manifold
$X\cup_{V}Y$ obtained by attaching smooth symplectic manifolds $X$
and $Y$ along a common symplectic submanifold $V$ of real codimension
two, as shown in Figure 3.  In this situation the generic fiber $Z_\lambda$ is a
symplectic smoothing of $Z_0$ and   is unique up   to symplectic deformation.  That deformation class is called the {\em symplectic sum} of
$(X,V)$ and $(Y, V)$.

 Alternatively, one can construct symplectic sums by hand.   Fix  a  symplectic
manifold $X$ with a symplectic  codimension two submanifold  $V$.
Suppose we are given a similar
pair $(Y,V)$ with a  symplectic  identification between the two
copies of $V$ such that,  under that identification, the normal
bundles
$N_X^V$  and $N_Y^V$ are dual as complex vector bundles.  We can then
glue $X$ to $Y$ along $V$ and ``round the corner'', essentially by
replacing the crossing $xy=0$ by $xy=\lambda$ in coordinates normal
to $V$.  For $\lambda\neq 0$ the  resulting manifold $Z_\lambda$ is
smooth and symplectic --- this is the symplectic sum.  In fact,  this
construction creates a family $\{Z_\lambda \}$  that together form a
symplectic fibration as above.

The  ``Symplectic Sum Formula'' developed in \cite{ip2} and
\cite{ip3}  expresses the GW invariants of  the symplectic sum
$Z_\lambda$ in terms of invariants of $X$ and $Y$.  The derivation
begins by considering what happens to $J$-holomorphic maps into
$Z_\lambda$ as $\lambda \to 0$.  By the compactness theorem for
$J$-holomorphic maps (for maps into  ${\mathcal Z}$), these limit to
maps into $Z_0=X\cup_{V}Y$, which can be separated into components on
the $X$ and $Y$ sides.   However, there are several complications.

    \begin{figure}[here]
    \begin{center}  \includegraphics[scale=.6]{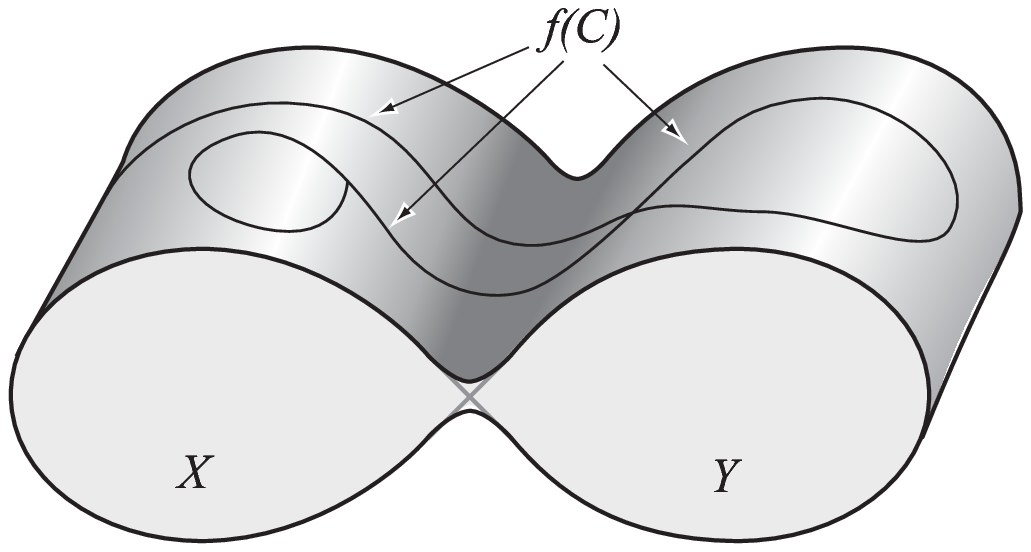}\end{center}
\caption{ Images $f(C)$ of $J$-holomorphic maps
   into $Z_\lambda$ as $\lambda\to 0$.}
   \end{figure}

    One complication is apparent in Figure 4:   connected curves in
$Z_\lambda$ may limit
to curves whose restrictions to $X$ and $Y$ are not connected, as
Figure 4 shows.  For
that reason the GW invariant, which counts stable curves from a
connected domain, is not the appropriate invariant for expressing a
sum formula.  Instead one should work with the ``Gromov-Taubes''
invariant $GT$, which counts stable maps from domains that need not be
connected.   This is not a substantive change:  one can move between
the GW invariants and the GT invariants contain the same information
and there are combinatorial formulas relating them.

A more substantial complication arises from the requirement, needed
for the compactness theorem mentioned above, that the
   the almost complex structures on $X$ and $Y$ match along $V$ and
extend smoothly to the fibers of the fibration ${\mathcal Z}$.    For
that reason, it is necessary to restrict to pairs  $(J,\nu)$ that
satisfy a certain
``$V$-compatibility'' condition (see \cite{ip3}).  That condition
implies, in particular, that $V$ is a $J$-holomorphic submanifold.
There is a price to pay for this specialization.  Because we can no
longer perturb $(J,\nu)$ along
$V$ at will, we do not have enough freedom to arrange that the limit curves are
transverse to $V$.  In fact, the images of  the
limit maps into $Z_0$ can meet $V$
at points with multiplicities and, worse, may have components
may be mapped entirely into $V$.

                     To count stable maps into $Z_0$ we look first at
maps into $X$ that have no
components  mapped into $V$.  These form a moduli space which is the
union of components
${\mathcal M}_s^V(X)$ labeled by the multiplicities
$s=(s_1,\dots, s_\ell)$ of the  intersection points with $V$.  Each
${\mathcal M}_s^V(X)$   can be compactified.  Those compact moduli spaces come
with evaluation maps like (\ref{1,evalmap}) that keep  track of the
domain curve, the images of the marked points, and the relative
homology
class of the image (details are given in \cite{ip2}).  Repeating the
construction of GW invariants, we then obtain  {\em Relative
Gromov-Witten invariants} $GW^V_X$ that count the number of stable
maps into $X$ that meet $V$ with specified multiplicities.

A third complication is the fact that the squeezing
process is not injective.  A precise analysis reveals that each map
into $Z_0$ that meets $V$ with
   multiplicities $s=(s_1,\dots, s_\ell)$ is the limit of $s_1\cdot
s_2\cdots s_\ell$ distinct maps into $Z_\lambda$.
    The maps within such a ``cluster''  have almost the same
image  but have different  domains.
As $\lambda\to 0$ the cluster coalesces, limiting to a  single map.

The limit process is reversed by a  standard
approach that goes back to Taubes and Donaldson: construct a space of
approximately
holomorphic maps (which incorporates the clustering phenomenon),
then correct the
approximate maps to true
holomorphic maps using  a fixed point theorem.   The upshot is a
precise statement, at the level of moduli spaces, of the fact that
stable maps into $Z_\lambda$ for small $\lambda$ are in one-to-one
correspondence to the clusters associated with pairs of  maps $(f,g)$
into $X$ and $Y$ that intersect $V$ at the same points with the same
multiplicities.

There remains the  issue of curves ``sinking into the neck'' in the
limit $\lambda\to 0$.   This is dealt with using a renormalization
argument.  For that, we consider the ruled manifold ${\Bbb P}_V={\Bbb P}({{ \Bbb C}}\oplus N_X^V)$, which is the ${\Bbb P}^1$ bundle over $V$ obtained by
adding an infinity section to the normal bundle to $V$.  We can then 
glue $X$ to $Y$ through a series of $n$ copies of ${\Bbb P}_V$, obtaining a  space which is singular along $n+1$ copies of $V$.  Introducing  smoothing parameters $\lambda_1,\dots,\lambda_{n+1}$ and applying the symplectic sum construction gives a fibration ${\mathcal Z}_n\to D^{n+1}$.   This  is similar to the fibration of Figure 3, but now the fibers near the central fiber are
manifolds $Z_\lambda$ with $n+1$ necks, as shown in Figure 5 for $n=2$.
 \begin{figure}[here]
  \centering
 \includegraphics[scale=.7]{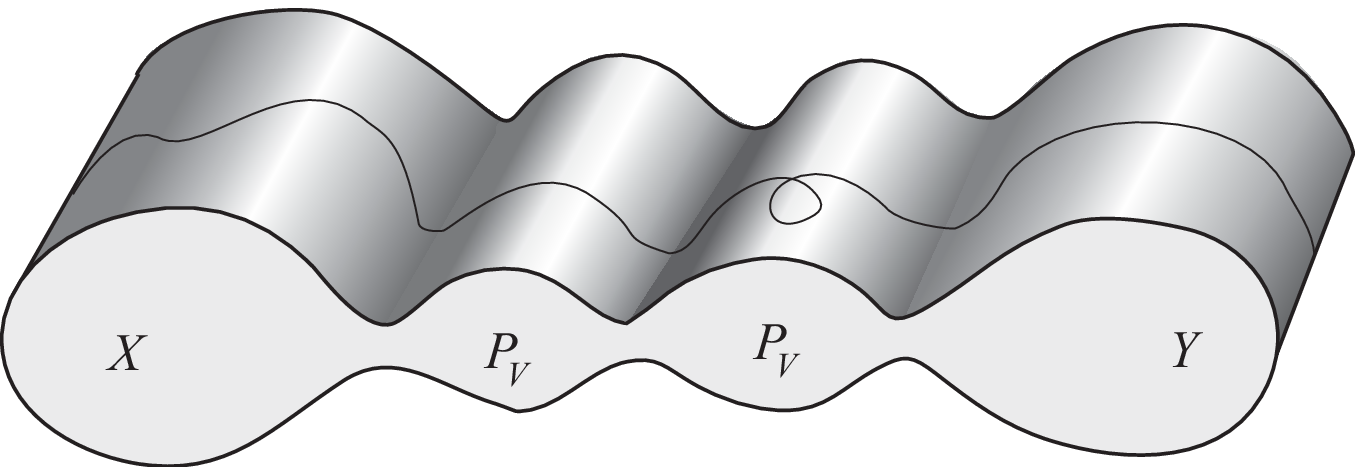}
\caption{Renormalizing along $V$ several times catches prevents curves from sinking into the neck.}
  \end{figure}
\noindent An energy bound shows that,  when one takes $n$
sufficiently large and lets $\lambda\to 0$, the limit maps have no
components sinking into most
necks of $Z_0$. Squeezing some or all of those necks decomposes the curves
in $Z_\lambda$ into curves in $X$ joined to curves in $Y$ by a chain of
curves in intermediate spaces ${\Bbb P}_V$, and the previous analysis
applies to each.  Working
through the combinatorics, one finds that
the total contribution of the entire region between $X$ and $Y$
is given by a certain $GT$ invariant $S_V$ of ${\Bbb P}_V$.  This
``$S$-matrix''  $S_V$  keeps track of how the genus, homology class,
and
intersection points with $V$ are ``scattered'' as the images of
stable maps pass
through the intermediate regions of Figure 5.

\bigskip\medskip

\noindent{\bf Symplectic Sum Theorem}{\em \quad   Let $Z$ be the symplectic
sum of $(X,V)$ and  $(Y,V)$.
Then the Gromov-Taubes invariant of $Z$ is given
in terms of the relative invariants of $X$ and $Y$  by
\begin{equation}
GT_Z\ =\   GT_X^{V} \, * \,
S_V\, * \, GT_Y^{V}
\label{3.SSF}
\end{equation}
where  $*$ is a certain `convolution operation'  and
   $S_V$ is a certain GW invariant of the ruled manifold ${\Bbb P}_V$.
}

\vskip 1cm

%%%%%%%%%%%%%%%%%%%%%%%%%%%%%%%%%%%%%%%%%%%%%%%%%%%%%%%%%%%%
%%%%%%%%%%%%%%%%%  Section 5  %%%%%%%%%%%%%%%%%%%%%%%%%%%%%%
%%%%%%%%%%%%%%%%%%%%%%%%%%%%%%%%%%%%%%%%%%%%%%%%%%%%%%%%%%%%%%

\setcounter{equation}{0}
\section{Application: Invariants of $E(1)$ }
\label{section5}
\bigskip

The formulas of the previous sections describe the behavior of GW invariants
under various cut-and-paste operations.
We referred to those as ``compuational tools'',
but have not yet done any computations or shown that
they have practical value.
We now remedy that  by presenting a striking application:
calculating GW invariants for the rational elliptic surface $E(1)$.

$E(1)$ is the simplest K\"{a}hler surface that fibers over the 2-sphere with
elliptic curves as fibers (that is, the fibers are tori with
varying complex structure).
   One standard construction goes as follows.
    Choose two generic degree 3 curves in ${\Bbb P}^2$.
     These are elliptic curves,
     given as the zero sets of cubic polynomials
     $f,g$ in homogeneous coordinates, which intersect at 9 points.
     Blowing up ${\Bbb P}^2$ at those 9 points gives the surface $E(1)$.
     It comes with 9 exceptional curves.
     Furthermore, the curves in the linear system
     $\{f+tg\, |\, t\in { \Bbb C}\}$ lift to a family of elliptic curves
     $\{E_t\}$  in $E(1)$ so that (a) $E_t$ and $E_{t'}$ are disjoint
     for $t\neq t'$, (b) every point of $E(1)$ lies on one and only $E_t$,
     and (c) each $E_t$ intersects each exceptional curve
     in one point with multiplicity one.
     Thus the $E_t$ are the fiber of a fibration $E(1)\to S^2$ and
     each exceptional curve is a section of that fibration.
     For generic cubic polynomials $f$ and $g$,
     this fibration has exactly 12 nodal fibers as shown in Figure 6 below.

Varying the blowup points gives different complex structures on $E(1)$ but,
up to deformation, there is a single symplectic structure on $E(1)$.
The homology $H_2(E(1);{\Bbb Z})$ is generated by the 9 section classes $S_i$
represented by the exceptional curves and by the class $F$ of a generic fiber.
These satisfy
$$
F\cdot F=0, \qquad F\cdot S_i=1, \qquad S_i\cdot S_j=-\delta_{ij},
$$
and altogether the intersection form of $E(1)$ is the sum $E_8\oplus H$ of
the $E_8$ matrix and the $2\times 2$ hyperbolic matrix.
The adjunction formula shows that the canonical class is $K=-F$.
The GW invariant  of the fiber class is  especially important: it is
$$
GW_{F,1}(E(1)) = 1
$$
(see Example 9.3 of \cite{ip2}).

\medskip

Now fix one section $S=S_i$ and consider the classes $S+dF$ for $d\in
{\Bbb Z}$.  For generic $(J, \nu)$, the space of stable maps $S^2\to E(1)$
form a moduli space  with dimension given by (\ref{1.dim}):
$$
\mbox{dim}\ {\overline{\mathcal M}}_{0,0}(E(1),S+dF)\ =\ 2\left[-1- 
K\cdot (S+dF)\right] \ =\ 0.
$$
Thus for each $d$ this moduli space is a finite set of points (with
sign).  We wish to compute the GW invariants that count  the signed
points in those moduli spaces.

\medskip

\noindent{\bf Problem A}\ \ Find the genus 0 Gromov-Witten
invariants $GW_{S+dF,0}$ of $E(1)$.   Equivalently, find the
generating function
\begin{equation*}
   F(t)\,=\,\sum_{d\geq 0} GW_{S+dF,0}\ t^{d}.
\end{equation*}

Notice that any class $S+dF$ with $GW_{S+dF,0}\ne 0$ is
represented by a $(J, \nu)$-holomorphic map for generic $(J,\nu)$.
As $\nu\to 0$ those maps limit to a $J$-holomorphic map
representing $S+dF$. By also varying $J$, we see that $S+dF$ has a
holomorphic representative for the $J$ of a complex elliptic
structure on $E(1)$, for which we know that $S$ also  has a
holomorphic representative.  Since distinct holomorphic curves
have non-negative intersection we conclude that $d=0$ or
$d-1=S\cdot(S+df)\geq 0$.  Thus there is no lose of generality  in
restricting the above sum to $d\geq 0$.

\medskip

There is a similar problem for genus 1 invariants.  Note that, by
(\ref{1.dim}),   the space of stable maps from a genus $g=1$ domain
into $E(1)$ has dimension 2.  We can get a GW invariant by adding a
marked point (thereby increasing the dimension of space of stable
maps to 4) and imposing the codimension 4 constraint  $\tau(F)$.  For that purpose, fix a smooth generic fiber $F_0\subset E(1)$. 

\medskip

\noindent{\bf Problem B}\ \ Find the genus 1 GW  invariants
$GW_{S+dF,1}(\tau(F))$ of $E(1)$. Equivalently, find the
generating function
\begin{equation*}
   H(t)\,=\,\sum_{d\geq 0} GW_{S+dF,1}(\tau(F)) \ t^{d}.
\end{equation*}

\bigskip

We will solve Problems A and B simultaneously by relating the
generating functions in two ways:  using a TRR formula and using the Symplectic
Sum formula.  That will give two differential equations involving
$F(t)$ and $H(t)$ that, together, determine
$F(t)$ and $H(t)$  explicitly.

First, we can expand the constraint $\tau(F)=\psi_1 \cup ev^*(F)$
using the TRR formula (\ref{trr}):
\begin{equation}
\label{5.picture}
GW_{S+dF,1}(\tau(F))   \ =\ \frac{1}{12}\
\begin{minipage}{.45in}
\includegraphics{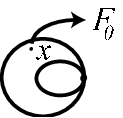}
\end{minipage}
    +\
\begin{minipage}{.47in}
\includegraphics{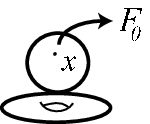}
\end{minipage}
\end{equation}
where the arrows to $F_0$ indicate that the marked points must be
mapped onto $F_0$.  The  pinched torus picture corresponds to a genus
0 domain with automorphism group ${\Bbb Z}_2$.  When that is expanded using
(\ref{splittingformula2}),  summed over $d$ and simplified using the
Divisor Axiom, we have
\begin{eqnarray}
\label{5.2}
\frac{1}{12}\cdot \frac12 \sum_{d\geq 0} (S+dF)^2 \
GW_{S+dF,0}(F)\ t^d
& = & \frac{1}{24} \sum_{d\geq 0} (2d-1)\ GW_{S+dF,0}(F)\ t^d\nonumber \\
& = &\frac{1}{24}\left(2t F'(t) - F(t)\right).
\end{eqnarray}
Maps corresponding to  the last picture of (\ref{5.picture}) have
components $f_1:S^2\to E(1)$ and $f_2:T^2\to E(1)$ that together
represent $S+dF$.  Composing with $\pi:E(1)\to S^2$ gives holomorphic
maps $\pi\circ f_i$ to $S^2$ whose total degree is 1. This implies
that one of the maps $\pi\circ f_i$ has degree  0, so its image  is a
single point.  Thus either $f_1$ or $f_2$ is a map into a single
fiber of $E(1)$.  But that cannot be the case for $f_1$ because
$f_1$, being a stable map with unstable domain, must represent a
non-trivial homology class and there are no homologically non-trivial
maps $S^2\to T^2$.  We conclude that $f_2$, if not trivial, is a
cover of a $J$-holomorphic fiber --- which must be $F_0$  ---  and
represents $d_2F$, while $f_1$ represents $S+d_1F$ with $d_1+d_2=d$.

Now applying Proposition 2.1 and the Divisor Axiom, we see that
the maps corresponding to  the last picture of
   (\ref{5.picture}) contribute
\begin{align}
\label{5.3}
   &\sum_{d\geq 0} t^d \sum_{d_{1}+d_{2}=d} \sum_{a}
    GW_{S+d_{1}F,0}(F,H_{a}) \cdot
    GW_{d_{2}F,1}(H^{a})  \nonumber \\\
   =&\ \sum_{d\geq 0} t^d\,\Big(
\sum_{k>0}
    GW_{S+(d-k)F, 0} \cdot
    k\,GW_{kF,1}\ + \
    \sum_{a}GW_{S+dF}(H_{a})\, GW_{0,1}(H^{a}) \Big)
   \nonumber \\\
   =&\   F(t)\ \sum_{d > 0}  d\,GW_{dF,1}\ t^{d} \,+\,
    \sum_{d\geq 0}t^{d}\,\sum_{a} GW_{S+dF}(H_{a})\,GW_{0,1}(H^{a}).
\end{align}
The invariant $GW_{0,1}(H^{a})$ is a count of degenerate maps and one
can show by hand that
\begin{equation}\label{pt-map}
   GW_{0,1}(H^{a})\,=\,\frac{1}{24}\,(H^{a}\cdot K)\,=\,
      -\frac{1}{24}\,(H^{a}\cdot F)
\end{equation}
(\,see \cite{ip3} \S 15\,). The invariant $d\,GW_{dF,1}$ counts
the $d$-fold coverings of a torus by tori and hence it is the sum
of all divisors of $d$, which is denoted by $\sigma(d)$. We set
\begin{equation}
  \label{G(t)}
   G(t)\ =\ \sum_{d=1}^\infty \sigma(d) \,t^d\ =\
   \sum_{k=1}^\infty \frac{kt^k}{1-t^k}
\end{equation}
Then, combining the last several equations yields a differential
equation that encodes the ``splitting the domain'' TRR relation:

\begin{equation}\label{F-trr}
    H(t)\,=\,\frac{1}{12}\left( \,t\,F^{\prime}(t)  - F(t)\right)\, +\,
G(t)\,F(t).
\end{equation}

\bigskip

We will next derive a similar
differential equation by splitting the target.   Let $V$ be one of
the fibers of $E(1)\to S^2$.   Setting $E(0)= S^2\times T^2$, let
$V'$ be a  fiber of $E(0)\to S^2$.  Since $V$ and $V'$ have trivial
normal bundles, we can identify them and form the symplectic sum.
That allows us to regard $E(1)$ as $E(1)\#E(0)$, as in Figure 6. We
will  apply the symplectic sum formula (\ref{3.SSF}), taking the
constraint $\tau(F)$ with $F$ represented by a fiber $F_0$ on the
$E(0)$ side.

   \vskip.4cm
   \hspace*{.6cm}
    \begin{center}
\begin{minipage}{2.5in}
   \includegraphics[scale=.9]{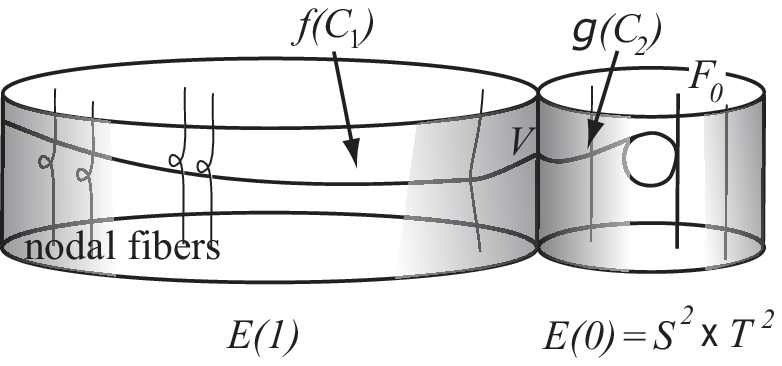}
   \end{minipage}
\hspace*{1.6cm}
\begin{minipage}{1.5in}
    \begin{tabular}{c}
  {\small{\bf Figure 6}}
  \addtocounter{figure}{1} \\ \\
Pinching  $E(1)\#_V E(0)$
    \end{tabular}
   \end{minipage}
   \end{center}
\vskip.2cm

In this case the symplectic sum formula simplifies considerably.  To
see how, one should think about the proof, rather than the formula
itself.  As we squeeze the neck in Figure 6 by letting $\lambda \to
0$,  maps from genus 1 domains that represent $S+dF$ limit to maps
$f:C_1\to E(1)$ and $g:C_2\to E(0)$  with (genus $C_1$)+ (genus $C_2$)=1
and whose images meet at a point of $V$ and with the image of $g$
contacting $F_0$ in a way that satisfies the constraint $\tau(F)$.

First consider the
possibility that $C_2$ has genus 0.  Then $g(C_2)$ represents $S+d'F$
for some $d'$, where $S$ and $F$ are the section and fiber classes of
$E(0)$.  Composing $g$ with the projection $E(0)\to T^2$ gives a map
$S^2\to T^2$, which must be zero in homology, so $d'=0$.
Stability implies that $C_2$ is a smooth 2-sphere
representing $S$.  Furthermore,  $F_0$ has trivial normal bundle and the image of $g$ has no components in $F_0$.  As we saw in Section 3, in that situation $\tau(F)$ is the constraint that the image of $g$ contacts $F_0$ to second order.  But $S\cdot F=1$, so that is impossible  for a fibered complex structure on
$E(0)$.  A  limiting argument shows that the same is true for generic
$J$, so  $GW_{S,0}(\tau(F))=0$.  We conclude that $C_1$ has genus 0
and $C_2$ has genus 1.

Exactly the same argument shows that the limit curve cannot have any
curves in the middle regions shown in Figure 5.  Consequently, the
matrix $S_V$ in the symplectic sum formula (\ref{3.SSF}) is the
identity and the relative and absolute invariants are the same.
Furthermore, since the maps intersect $V$ at a single point with
multiplicity 1 the domains $C_1$ and $C_2$ are connected and  the
multiplicity issues mentioned in Section 4 do not arise, and the
convolution in (\ref{3.SSF})  is simply the product of the generating
functions. Since $C_1$ has genus 0 the generating function on the $E(1)$
  side is just $F(t)$. Therefore, we have
\begin{equation}\label{SumF}
  H(t)\ =\ F(t)\,\sum_{d\geq 0}\,GW_{S+dF,1}(E(0))(\tau(F))\,t^{d}
\end{equation}

We can calculate the generating function for $E(0)$ by applying  TRR
exactly as in (\ref{5.2}) and (\ref{5.3}).  In fact, since $S^{2}=0$
the contribution of the first picture of
(\ref{5.picture}), given by the left-hand side of (\ref{5.2}), is
zero.  Moreover,   the invariant $GW_{S+dF,0}(E(0))=0$ for $d>0$, so
the  contribution of the second picture of (\ref{5.picture}), is
given by the last line of (\ref{5.3}) with $F(t)$ replaced by 
$GW_{S,0}(E(0))=1$.
That is evaluated using (\ref{pt-map}), noting that the canonical
class of $E(0)$ is $-2F$, and again recognizing that $dGW_{dF,1}(E(0))$ is
twice the number of $d$-fold covers of a torus by a torus (since, for
generic $J$, $E(0)$ has exactly two holomorphic fibers in the fiber
class).
Consequently,
$$
\sum_{d\geq 0} GW_{S+dF,1}(E(0))(\tau(F))\ t^d  \
=\
-\frac{1}{12} + 2G(t).
$$

With that, the symplectic sum formula
corresponding to the decomposition of Figure 6  has the simple form

\begin{equation}\label{F-sum}
   H(t)\,=\,-\frac{1}{12}\,F(t)\  +\  2\,G(t)\,F(t).
\end{equation}

\medskip

Now the punch line:  equating formulas (\ref{F-trr}) and
(\ref{F-sum}) gives the ODE
\begin{equation*}
t\,F^{\prime}(t)\,  =\, 12\,G(t)\,F(t)
\end{equation*}
with the initial condition $F(t)=GW_{S,0}=1$. The
solution of this ODE is given by
\begin{equation}\label{E:Main}
   F(t)\,=\,\prod_{d\geq 1}
   \left(\frac{1}{1-t^{d}}\right)^{12}
\end{equation}
and we can then obtain $H(t)$ from (\ref{F-sum}).  With one more 
application of the sum formula we can obtain the corresponding 
generating function for curves of each genus $g>0$ in $E(1)$ (see 
\cite{ip1}).

\vskip 1cm

%%%%%%%%%%%%%%%%%%%%%%%%%%%%%%%%%%%%%%%%%%%%%%%%%%%%%%%%%%%%
%%%%%%%%%%%%%%%%%  Section 6  %%%%%%%%%%%%%%%%%%%%%%%%%%%%%%
%%%%%%%%%%%%%%%%%%%%%%%%%%%%%%%%%%%%%%%%%%%%%%%%%%%%%%%%%%%%%%

\setcounter{equation}{0}
\section{Enumerative Conjectures and the  Inadequacy of  GW Invariants}
\label{section6}
\bigskip

  The calculation done in Section 5 is very encouraging:  we 
successfully used symplectic gluing methods determine GW invariants 
in an interesting case, and the result automatically appeared in the 
elegant form (\ref{E:Main}).  Moreover, the result correctly 
enumerates curves in a case --- rational curves in the classes $S+dF$ 
in $E(1)$ ---  of interest to algebraic geometers. Obviously, one 
would like to extend this approach and develop it into a general 
method for calculating GW invariants for whole classes of symplectic 
manifolds.  One can further hope to relate those GW invariants to 
enumerative invariants of    K\"{a}hler surfaces.
 
  This idea is 
especially appealing because it seems to be linked to some remarkable 
conjectures in enumerative algebraic geometry (see \cite{go}).  These 
arose from the observation that in the cases where  enumerative 
invariants have been proved or conjectured (some by algebraic 
geometers and some by string theorists),  the generating functions 
are remarkably similar.  They are always products of a few special 
functions.  Two of those special functions, namely $F(t)$ in 
(\ref{E:Main}) and  $G(t)$ in (\ref{G(t)}), arose naturally in gluing 
computations of the previous section.  Thus it is reasonable to 
conjecture that there are similar  universal formulas for the GW 
invariants of  K\"{a}hler surfaces.
  
  But there is a catch.  This 
vision of simple, calculable and non-trivial GW generating functions 
is not true for the usual Gromov-Witten invariants.  In fact, for 
K\"{a}hler surfaces with geometric genus $p_g>0$ {\em all of the GW 
invariants of interest  vanish}.   This phenomenon has a simple 
explanation, which we describe in this section.  The discussion here 
is preparation for introducing the ``family GW invariants'' in the 
next section.   
  
  Recall that the geometric genus $p_g$ of a 
K\"{a}hler surface is the (complex) dimension of the space 
$H^{0,2}(X)$ of  holomorphic (0,2)-forms.  The role of $p_g$ can 
already be seen in the following example, which involves the elliptic 
surface $E(n)$ obtained  as the $n$-fold fiber sum of copies of 
$E(1)$ as shown in Figure 7. 
  
   \begin{figure}[here]
\centering
\includegraphics[scale=.8]{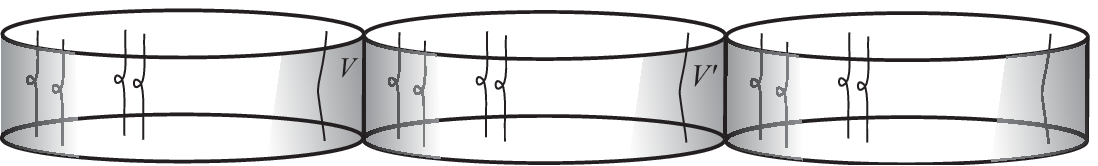}
\caption{$E(3)=E(1)\#_V E(1)\#_{V'} E(1)$}
\end{figure}

 \begin{example}
Let  $E(n)\to S^{2}$ be an  elliptic surface with a
section of self-intersection number $-n$.  This is a simply-connected 
K\"{a}hler manifold with a fiber class $F$  and  a section class $S$ satisfying 
 $S\cdot S=-n$. Using the adjunction formula 
and Riemann-Roch, one sees that  the canonical class   is $K=(n-2)F$ 
and the   geometric genus is  $p_g=(n-1)$.  Then by (\ref{1.dim}) the 
dimension of the GW moduli space for
the class $S+dF$ with genus $g=0$ and no marked points is
\begin{equation*}
\mbox{dim}\ {\overline{\mathcal M}}_{g,0}(E(n),S+dF)\ =\ 
2(\,-K\cdot(S+dF)-1)\,=\,2(1-n)=-2p_{g}.
\end{equation*}
This is negative when $n\geq 2$, meaning that the moduli space is 
empty for generic $J$.  Thus  the genus $0$ GW-invariants for the 
class $S+dF$ are all zero.
 \end{example}

\medskip
  
 More generally, in a K\"{a}hler surface $X$, the space 
of holomorphic curves representing a fixed homology class $A$ forms a 
space called the   (generalied) Severi variety (see \cite{bl3}). 
When the Severi variety  has dimension $n$, one can impose $n$  point 
constraints by introducing marked points $x_1, \dots x_n$ on the 
domain and $n$ points $y_1, \dots y_n$ in $X$ and restricting 
attention to the set of maps $f$ with $f(x_j)=y_j$ for $j=1,\dots n$. 
When certain tranversality conditions are met, that will be a finite 
set of maps, and count of those maps is an enumerative invariant of 
$X$.

  However, when $p_g>0$ that enumerative invariant is not 
a GW invariant.  The reason is that the (formal)
dimension of the GW moduli space is strictly less than the
dimension of the Severi variety.  Here is a specific computation.

\begin{example}
\label{ex6.1}
{\rm Let $X$ be a
K\"{a}hler surface with canonical class $K$.  Suppose that
$f:C\to X$
is an embedded holomorphic map representing a class $A$ with $A-K$
ample.
Then the space ${\mathcal V}_A$ of curves representing $A$ has
complex dimension
$$
\mbox{dim}_{ \Bbb C} \ {\mathcal V}_A\ =\ d_A
+p_g
$$
where $d_A$ is the dimension (\ref{1.dim}) of the GW moduli
space of maps representing $A$.}
\end{example}

\begin{proof} Let
$\xi$ the line bundle  over $X$
corresponding to the divisor $f(C)$ and let $N$ be the normal bundle.  Fixing a section with divisor $f(C)$ gives a 
sheaf exact sequence
$0\,\to\,{\mathcal O}_{X}\,\to\,\xi\,\to\,N\,\to\,0
$, which induces  the  long exact sequence
\begin{equation*}
   \cdots\,\to\,H^{1}(X;\xi)\,\to\,H^{1}(C;N)\,\to\,
   H^{2}(X;{\mathcal O}_{X})\,\to\,H^{2}(X;\xi)\,\to\,0.
\end{equation*}
But since $A-K$ is ample,  it  follows from Serre duality and
Kodaira Vanishing Theorem that  $H^{1}(X;\xi)$ and
$H^{2}(X;\xi)$ both  vanish.  Thus
$H^{1}(C;N)\cong  H^{0,2}(X)$ has complex dimension $p_g$.  On the other
hand, the tangent space to ${\mathcal V}_A$  is the space
$H^0(C,N)$ of holomorphic deformations of $f(C)$.  Finally,
by
Riemann-Roch and (\ref{1.dim})
$$
d_A\ =\ H^0(C,N) - H^1(C,N)\
=\ \dim_{ \Bbb C} {\mathcal V}_A -p_g.
$$
\end{proof}

The situation in
Example \ref{ex6.1} is built into the algebraic geometry viewpoint of 
Gromov-Witten invariants. Ignoring technical difficulties, the 
picture is as follows.   On a  K\"{a}hler surface, the space of 
stable maps ${\overline{\mathcal M}}_{g,k}(X,A)$  (which is the 
Severi variety when the maps are embeddings) is a compact projective 
variety of some dimension $n$.   There is a rank $p_g$ bundle 
${\mathcal{E}}$ over that space 
  \begin{equation}
  \begin{CD}
{\mathcal{E}}\\
  @VVV\\
  {\overline{\mathcal M}}_{g,k}(X,A)
\end{CD}
  \label{6.E}
\end{equation}
  whose fiber over a map $f$ is 
$H^{0,1}(C;f^*N)$.  The Euler class of ${\mathcal{E}}$ defines a 
homology class
\begin{equation}
{\overline{\mathcal 
M}}^{virt}_{g,k}(X,A)\overset{def}{=} [{\overline{\mathcal 
M}}_{g,k}(X,A)]\cap e({\mathcal{E}})  \ \in\  H_{2(n-p_g)}({\mathcal 
V}_A)
\label{6.VMC}
\end{equation}
 called the virtual moduli cycle. 
Algebraic geometers define GW invariants using this in place of the 
space (\ref{1.Mstable}) of perturbed holomorphic maps.

 To 
relate this picture to  the geometric analysis viewpoint, consider a 
perturbation $\nu$ as in (\ref{1.pertholomorphicmapeq}).  At each map 
$f$, the  linearization   of the $J$-holomorphic map equation is an 
operator $L_f:\Omega^0(f^*TX)\oplus H^{0,1}(TC) \to 
\Omega^{0,1}(f^*TX)$ (see \cite{rt2}).  The projection of $\nu$ into 
the cokernel  $L_f$  defines a section of ${\mathcal{E}}\to 
{\overline{\mathcal M}}_{g,k}(X,A)$.  The zeros of that section are 
exactly the $J$-holomorphic maps $f$ that, by the Implicit Function 
Theorem, can be uniquely corrected to solutions of the perturbed 
holomorphic map equation $\overline{\partial}_Jf=t \nu$ for small 
fixed $t$.  That zero set thus represents the virtual 
  moduli cycle 
(\ref{6.VMC}) and is homologous to the moduli space (\ref{1.Mstable}) 
defined by the perturbed equations.

  The distinction between 
the GW and the enumerative counts should now be clear.  Enumerative 
geometers impose constraints and count curves.  That is in principle 
the same as a count of  constrained maps in the space 
${\overline{\mathcal M}}_{g,k}(X,A)$ of stable maps  for a K\"{a}hler 
$J$.  But it is entirely different from the GW invariants, which 
count constrained maps in the codimension $p_g$ subspace 
${\overline{\mathcal M}}^{virt}_{g,k}(X,A)$, or equivalently in the 
space of perturbed holomorphic maps for generic $(J,\nu)$.

\medskip
  
This brings us to the key point: to  get at the 
enumerative invariants we need a new version of GW invariants.  Those 
invariants should be based on   a moduli space  of $J$-holomorphic 
maps that 
\medskip
\begin{itemize}
\item  is equal to the space of stable maps for a 
K\"{a}hler $J$, and \medskip
\item has  dimension $2p_g$ greater than the dimension 
of the usual  GW moduli space.
\end{itemize}
\medskip
There are several approaches to doing that.
Bryan 
and Leung \cite{bl1},\cite{bl2}
defined and calculated such invariants for K3 and Abelian surfaces by using the
twistor family.  Behrend-Fantechi \cite{bf} have defined
invariants for a more general class of algebraic surfaces using
algebraic geometry.   Here we will focus on the analytic approach of 
the first author.  That is described in the next section.

\vskip 1cm

%%%%%%%%%%%%%%%%%%%%%%%%%%%%%%%%%%%%%%%%%%%%%%%%%%%%%%%%%%%%
%%%%%%%%%%%%%%%%%  Section 7  %%%%%%%%%%%%%%%%%%%%%%%%%%%%%%
%%%%%%%%%%%%%%%%%%%%%%%%%%%%%%%%%%%%%%%%%%%%%%%%%%%%%%%%%%%%%%

\setcounter{equation}{0}
\section{$J_\alpha$-holomorphic Maps on K\"ahler surfaces}
\label{section7}
\bigskip

This section describes the first author's  family GW
invariants for K\"{a}hler  surfaces. These are based on the 
observation that, associated to each K\"{a}hler structure $J$, there 
is a  2$p_{g}$-dimensional family of almost complex structures, very 
naturally parameterized by $H^{0,2}(X)$.  While this parameter space 
is non-compact, one can show that in many interesting cases the 
resulting family moduli space is compact, and hence defines 
invariants.
Details of the construction outlined here can be found 
in \cite{l1}.

 Fix a compact K\"ahler surface
$(X,\omega,J, g)$ and choose the $2p_{g}$-dimensional parameter
space
\begin{equation*}\label{para}
   {\mathcal H}\,=\,\mbox{Re}\big(\,H^{2,0}\oplus  H^{0,2}\,\big).
\end{equation*}
Using the metric,  each $\alpha\in {\mathcal H}$ defines an endormorphism
$K_{\alpha}$ of $TX$ by the equation
   $$\langle u,K_{\alpha}v\rangle=\alpha(u,v).
$$
It follows that $K_{\alpha}$ is  skew-symmetric and anti-commutes
with $J$. For each map $f:C\to X$, we then have
   a 2$p_g$ dimensional family of  elements $K_\alpha\partial f j$
   in $\Omega^{0,1} (f^*TX)$.  We  modify the
$J$-holomorphic map equation (1) by considering the pairs
$(f,\alpha)$ satisfying
\begin{equation}
\overline{\partial}_Jf= K_\alpha\partial f j.
  \label{junhoeq}
\end{equation}

We can recast this as a family of pseudo-holomorphic maps
parameterized by ${\mathcal H}$. The above definitions imply that
$Id+JK_{\alpha}$ is injective, and hence invertible, for each $\alpha
\in {{\mathcal H}}$. Thus
\begin{equation*}
J_{\alpha} = (Id + JK_{\alpha})^{-1}J\,(Id + JK_{\alpha})
\end{equation*}
is a family of almost complex structures on $X$ parameterized by 
$\alpha$ in the 2$p_{g}$-dimensional linear space ${{\mathcal H}}$. 
Some straightforward algebra shows that (\ref{junhoeq}) is equivalent 
to
\begin{equation}
\overline{\partial}_{J_\alpha} f= 0.
  \label{junhoeq2}
\end{equation}

\medskip

Thus for each $J$ we get a family of holomorphic maps consisting of 
the pairs $(f,\alpha)$ where
   $f:C\to X$ is a map from a stable curve  to $X$ and 
$\alpha\in{{\mathcal H}}$ which satisfy (\ref{junhoeq}) or 
equivalently (\ref{junhoeq2}).  Fixing the genus $g$ and number $k$ 
of marked points on the domain and a class $A\in H_2(X; {\Bbb Z})$, 
the  set of such pairs forms a moduli space
   \begin{equation}
   {\overline{\mathcal M}}^{{{\mathcal H}}}_{g,n}(X,A,J) = \{\ (f,\alpha)\ |\
{\overline\partial_{J_\alpha}} f = 0\,,\
   \alpha\in{{\mathcal H}}\,,\ [f]=A  \}
    \label{modulispace}
\end{equation}
whose formal dimension
  \begin{equation}
\mbox{dim}\ {\overline{\mathcal M}}^{{{\mathcal 
H}}}_{g,k}(X,A)=2\left[(g-1)- K\cdot A+k+p_g\right]
\label{7.dim}
\end{equation}
  is $2p_g$ larger than the dimension (\ref{1.dim}) of the usual GW 
moduli space.  As in (\cite{rt1}),  we can add a perturbation term to 
the equation (\ref{junhoeq}) to make the moduli space an orbifold of 
this dimension.  Then, following the standard procedure described in 
Section 1, we can  pushforward by the evaluation and obtain 
symplectic invariants {\em provided the modui space is compact}.

   Compactness is an
issue because the
parameter $\alpha$, which  ranges over the vector space ${{\mathcal 
H}}$,  may be unbounded.  Nevertheless, one can use specific 
geometric arguments to
show that the moduli space {\em is} compact for many classes $A$,
and consequently the invariants are well-defined  for those
classes. In \cite{l1},  we proved that is true for many
cases important in algebraic geometry.

\begin{theorem}[\cite{l1}]
\label{7.existence}
  The family GW invariants $GW^{{{\mathcal H}}}_{A,g}$ are well-defined
\begin{enumerate}
\item[(a)]  for all non-zero classes $A$ on a K3 or Abelian surface.
\item[(b)] for all (1,1) classes $A$ on  a minimal  elliptic surface of
   Kodaira dimension 1 satisfying   $A\cdot (\mbox{fiber class})>0$, and
\item[(c)] for many $A$ on a  K\"ahler surface of general type
   (specifically, for   (1,1) classes $A$ that are not  linear combinations of
   components of the
   canonical class $K$).
\end{enumerate}
\end{theorem}

For K3 and abelian surfaces, our family GW-invariants are identical to the
``twistor family invariants''  used by Bryan and Leung. In particular,
(i) the invariants are independent of complex structures, (ii)
they count actual holomorphic curves representing primitive
classes, and (iii) for any two primitive classes $A$ and $B$  with 
$A^2=B^2$   there
is a orientation preserving diffeomorphism $h:X\to X$ such that
$h_{*}A=B$ and
\begin{equation}\label{BL}
   GW^{{{\mathcal 
H}}}_{mA,g}\big(h^{*}(\gamma_{1}),\cdots,h^{*}(\gamma_{k})\big)\,=\,
   GW^{{{\mathcal H}}}_{mB,g}\big(\gamma_{1},\cdots,\gamma_{k}\big).
\end{equation}

\vskip 1cm

%%%%%%%%%%%%%%%%%%%%%%%%%%%%%%%%%%%%%%%%%%%%%%%%%%%%%%%%%%%%
%%%%%%%%%%%%%%%%%  Section 8  %%%%%%%%%%%%%%%%%%%%%%%%%%%%%%
%%%%%%%%%%%%%%%%%%%%%%%%%%%%%%%%%%%%%%%%%%%%%%%%%%%%%%%%%%%%%%

\setcounter{equation}{0}
\section{Splitting Formulas for Family Invariants}
\label{section8}
\bigskip

  One can
extend the splitting formulas of Section 2 to the  family GW invariants by 
exploiting the fact that the  
  parameter space ${{\mathcal H}}$ is 
a linear space.
For simplicity, we will describe the extension of the splitting
formula only for the case of elliptic surface $E(n)$ with the
classes $S+dF$, $d\in {\Bbb Z}$.  The general case is done in  \cite{l1}.

  \smallskip

Let $\sigma$ be the gluing map of Figure 1.  Each 
$J_{\alpha}$-holomorphic map whose domain lies in the   image of $\sigma$ 
can be written as a pair of $J_{\alpha}$-holomorphic maps 
$(f_{1},\alpha):C_{1}\to X$ and $(f_{2},\alpha):C_{2}\to
X$ with $f_1(x_{k_{1}+1})=f_2(y_{1})$.   Now it is a basic fact, 
proved in \cite{l1}, that any solution of (\ref{junhoeq}) 
representing a (1,1) class has $\alpha=0$.  Since $S+dF$ is a $(1,1)$ 
class
with respect to the complex structure of $E(n)$,  both $f_{1}$ and $f_{2}$ are
holomorphic. Thus, either $[f_{1}]=S+d_{1}F$ and
$[f_{2}]=d_{2}F$ or $[f_{1}]=d_{1}F$ and $[f_{2}]=S+d_{2}F$ with
$d_{1}+d_{2}=d$. Define
\begin{align*}
    p_{1,d_{1},d_{2}}\,:\,
    {\overline{\mathcal M}}_{g_{1},k_{1}+1}^{{{\mathcal 
H}}}\big(E(n),S+d_{1}F\big)\times
    {\overline{\mathcal M}}_{g_{2},k_{2}+1}^{{{\mathcal 
H}}}\big(E(n),d_{2}F\big)\,
    &\to\,{{\mathcal H}}\times{{\mathcal H}} \\
    p_{2,d_{1},d_{2}}\,:\,
    {\overline{\mathcal M}}_{g_{1},k_{1}+1}^{{{\mathcal 
H}}}\big(E(n),d_{1}F\big)\times
    {\overline{\mathcal M}}_{g_{2},k_{2}+1}^{{{\mathcal 
H}}}\big(E(n),S+d_{2}F\big)\,
    &\to\,{{\mathcal H}}\times{{\mathcal H}}
\end{align*}
by $p_{i,d_{1},d_{2}}
\left((f_{1},\alpha_{1}),(f_{2},\alpha_{2})\right)=(\alpha_{1},\alpha_{2})$
for $i=1,2$.
Then the moduli space of $J_{\alpha}$-holomorphic maps $(f,\alpha)$ from
nodal domains lying in $\mbox{Im}(\sigma)$ defines a
homology class
\begin{equation}\label{com1-1}
    \sum_{i=1,2}\sum_{d=d_{1}+d_{2}} \left(
    \big[\,
    p^{-1}_{i,d_{1},d_{2}}(\Delta^{{{\mathcal H}}})\, \big]\,\cap\, \Delta^{*}
    \right)
\end{equation}
where $\Delta^{{{\mathcal H}}}$ is the diagonal in ${{\mathcal H}}\times
{{\mathcal H}}$ and $\Delta^{*}=ev^*_{k_1+1,1}(PD[\Delta])$ as in
(\ref{compostion1}).  As in (\ref{GWCompostion1}), the family GW 
invariants with
constraint $PD(\sigma)$
\begin{equation*}
    GW^{{{\mathcal H}}}_{S+dF,g}(PD(\sigma);\gamma_{1},\cdots,\gamma_{k})
\end{equation*}
are obtained by evaluating (\ref{com1-1}) on the product of the 
classes $\gamma_{1},\cdots,\gamma_{k}$.

The  linearity of the parameter space ${{\mathcal H}}$
allows us to deform the spaces representing the summands of 
(\ref{com1-1}).  Set 
$\Delta_{t}^{{{\mathcal H}}}\, =\,
    \{\, (\,\alpha,t\alpha\,)\, |\, \alpha\,\in\,{{\mathcal H}}\, \}$ 
and consider the
spaces
\begin{equation}\label{com1-2}
    p_{1,d_{1},d_{2}}^{-1}(\Delta^{{{\mathcal H}}}_{t})\,\cap\,
    ev^{-1}_{k_1+1,1}(\Delta)
\end{equation}
for  $0\leq t\leq 1$.  These consist of triples 
$(f_1,f_2, \alpha)$ where $f_1$ is $J_\alpha$-holomorphic and $f_2$ 
is  $J_{t\alpha}$-holomorphic for the same $\alpha$.  By Theorem 7.1b 
the family moduli space
for the class $S+d_{1}F$ is compact.  That may not be true for the 
class $d_2F$, but nevertheless the spaces (\ref{com1-2}) are compact: 
the compactness of the space of $(f_1,\alpha)$ gives a bound on 
$|\alpha|$, which then implies the compactness of set of $(f_2,t 
\alpha)$.  Consequently, as $t$ varies the spaces (\ref{com1-2}) 
trace out a compact cobordism, so represent the same homology class. 
Thus we can replace the summands of (\ref{com1-1}) by the classes of 
(\ref{com1-2}) with $t=0$, namely
\begin{equation*}
    \big[
    {\overline{\mathcal M}}_{g_{1},k_{1}+1}^{{{\mathcal 
H}}}\big(E(n),S+d_{1}F\big) \big]
    \otimes  \big[
    {\overline{\mathcal M}}_{g_{2},k_{2}+1}\big(E(n),d_{2}F\big) \big]
    \,\cap \Delta^{*}
\end{equation*}
where the second factor is the ordinary GW moduli space.
The same reasoning applies with $p_1$ replaced by $p_2$, so 
(\ref{com1-1}) becomes
\begin{align*}
    &\sum_{i}\sum_{d=d_{1}+d_{2}} \left(
    \big[\,
    p^{-1}_{i,d_{1},d_{2}}(\Delta^{{{\mathcal H}}})\, \big]\,\cap\, \Delta^{*}
    \right) \notag \\
    %%%%%%%%%%%%%%%%%%%%%%%
    &=\,\sum_{d=d_{1}+d_{2}}\left(
    \big[
    {\overline{\mathcal M}}_{g_{1},k_{1}+1}^{{{\mathcal 
H}}}\big(E(n),S+d_{1}F\big) \big]
    \otimes  \big[
    {\overline{\mathcal M}}_{g_{2},k_{2}+1}\big(E(n),d_{2}F\big) \big]
    \,\cap \Delta^{*}\right) \notag \\
    &+\,\sum_{d=d_{1}+d_{2}}\left(
    \big[
    {\overline{\mathcal M}}_{g_{1},k_{1}+1}\big(E(n),d_{1}F\big) \big]
    \otimes  \big[
    {\overline{\mathcal M}}_{g_{2},k_{2}+1}^{{{\mathcal 
H}}}\big(E(n),S+d_{2}F\big) \big]
    \,\cap \Delta^{*}\right).
\end{align*}
We can then split  the diagonal  as in (\ref{splitdiag}) to obtain a 
splitting formula  for our  family
GW invariants.

\begin{proposition}
\label{8.cl} Let $\sigma$ be the gluing map of Figure 1 and fix 
cohomology bases  $\{H_{a}\}$ and $\{H^{a}\}$  as in 
(\ref{splitdiag}).    Then for any $\gamma _{1},\cdots ,\gamma_k\in 
H^*(E(n);{\Bbb Q})$
\begin{align*}
    & GW_{S+dF,g}^{{{\mathcal H}}}(PD(\sigma );\gamma _{1},\cdots 
,\gamma _{k}) \\
    & =\ \sum_{a}\sum_{d=d_{1}+d_{2}}
      GW_{S+d_{1}F,g_{1}}^{{{\mathcal H}}}(\gamma_{1},\cdots 
,\gamma_{k_{1}},H_{a})\,
      GW_{d_{2}F,g_{2}}(H^{a},\gamma_{k_{1}+1},\cdots ,\gamma_{k}) \\
    & +\ \sum_{a}\sum_{d=d_{1}+d_{2}}
      GW_{d_{1}F,g_{1}}(\gamma_{1},\cdots ,\gamma_{k_{1}},H_{a})\,
      GW_{S+d_{2}F,g_{2}}^{{{\mathcal 
H}}}(H^{a},\gamma_{k_{1}+1},\cdots ,\gamma_{k})
\end{align*}
where $GW$ denotes the ordinary GW invariants of $E(n)$.
\end{proposition}

\medskip

Note that in this formula we multiply family GW invariants by ordinary GW invariants.  Thus Proposition \ref{8.cl} can be used to ``bootstrap'' information about ordinary GW invariants in obtain information about family invariants.

\bigskip

The analogous formula for a non-separating node (as in Figure 2) is especially 
straightforward.   It is exactly the same as the corresponding formula (\ref{splittingformula2})
 for the GW invariant:
\begin{align}\label{com2}
     GW_{S+dF,g}^{{{\mathcal H}}}(PD(\theta );\gamma _{1},\cdots ,\gamma_{k})\,
     &=\,\sum_{a}GW_{S+dF,g-1}^{{{\mathcal H}}}(\gamma_{1},\cdots 
,\gamma_{k},H_{a},H^{a})
     \notag \\
     &=\,(S+dF)^{2}\,
     GW_{S+dF,g-1}^{{{\mathcal H}}}(\gamma_{1},\cdots ,\gamma_{k})
\end{align}
where the second equality follows by  the Divisor Axiom   (\ref{divisoraxiom}).

\vskip1cm

%%%%%%%%%%%%%%%%%%%%%%%%%%%%%%%%%%%%%%%%%%%%%%%%%%%%%%%%%%%%%%%%%%%%%%%%%%%%%
%%%%%%%%%%%%%%%%%%%%%% Section 9  %%%%%%%%%%%%%%%%%%%%%%%%%%%%%%%%%%%%%%%
%%%%%%%%%%%%%%%%%%%%%%%%%%%%%%%%%%%%%%%%%%%%%%%%%%%%%%%%%%%%%%%%%%%%%%%%%%%%%

\setcounter{equation}{0}
\section{Family Invariants for $E(n)$ }
\label{section9}
\bigskip

In this section, we will compute  family invariants
of $E(n)$ using the approach we took in Section 5  for  $E(1)$ --
combining  TRR and symplectic sum formulas.
We have just derived the relevant TRR formulas for the family invariants, so 
the next step is to similarly extend the Symplectic Sum Formula. 
Here we will do that for the relevant case of
elliptic surfaces.

The relevant sum formula is easily stated.  Fix a smooth fiber $V$ of 
an elliptic fibration  $E(n)\to S^{2}$ and, as before, regard $E(n)$ 
as the symplectic sum  $E(n)\#_{V} E(0)$. Choose a smooth bump 
function $\beta$ that vanishes in a small $\delta$-neighborhood of 
$V$ and is 1 everywhere outside a $2\delta$-neighborhood of $V$.  We 
can define a set of
relative family GW invariants for $\big( E(n),V\big) $ by the
simple method of replacing $\alpha\in {{\mathcal H}}$  by 
$\beta\alpha$ in (\ref{junhoeq}). Then each $J_{\beta\alpha}$ 
satisfies $J_{\beta\alpha}=J$ in a
neighborhood of $V$. Using simple properties of $J_\alpha$-holomorphic maps, one
can show that the resulting relative family moduli space for the 
classes $S+dF$ with $d\in {\Bbb Z}$ is
compact.
It then follows from the Symplectic Sum Theorem in Section 4 and the 
exactly same arguments of Section 5 that
\begin{equation}\label{FVS}
  GW_{S+dF,1}^{\mathcal H}(\tau(F))\ =\ 
\sum_{d_{1}+d_{2}=d}\,GW_{S+d_{1}F,0}^{\mathcal H}\,\cdot \,
  GW_{S+d_{2}F,1}(\tau(F)).
\end{equation}
In fact, the proof of the sum formula given in \cite{ip3} holds, 
essentially without
change,  for the family invariants provided that all the maps that 
arise in the limit
$Z_\lambda \to Z_0$ lie in {\em compact}  relative moduli spaces.
That is true for a fibration $\lambda : {\mathcal Z} \to D$  whose generic
fiber is $E(n)$ and whose center fiber is $E(n)\cup_V E(0)$ because 
the relative family moduli space is compact.

\bigskip

The splitting formulas in Section 8 and the sum formula (\ref{FVS})
allow us to use the arguments of Section 5 to compute the
genus
$0$ family GW invariants $GW^{{{\mathcal H}}}_{S+dF,0}$ of $E(n)$ for 
the class $S+dF$.  For these invariants, the family moduli space 
(with no marked points)
has dimension $2(1-n)+2p_{g}=0$. As before, we assemble the
invariants in the generating function
\begin{equation*}
    F(t)\,=\,\sum_{d\geq 0}
    GW^{{{\mathcal H}}}_{S+dF,0}\,t^{d}
\end{equation*}
  and compare it with the generating function for the genus $g=1$
family invariants with the descendent constraint $\tau(F)$
\begin{equation*}
    H(t)\,=\,\sum_{d\geq 0}
    GW^{{{\mathcal H}}}_{S+dF,1}\big(\tau(F)\big)\,t^{d}.
\end{equation*}

Applying  the TRR formula as in (\ref{5.picture})  and (\ref{5.2}) and 
simplifying using
the splitting formulas in Section 8    gives
\begin{align*}
   H(t)\,&=\,\frac{1}{12}\cdot \frac12 \sum_{d\geq 0} (S+dF)^2 \
             GW_{S+dF,0}^{{{\mathcal H}}}(F)\ t^d \\
         &+\,\sum_{d\geq 0}\,\Big(\,\sum_{d_{1}+d_{2}=d} \sum_{a}
             GW_{S+d_{1}F,0}^{{{\mathcal H}}}(F,H_{a}) \cdot
             GW_{d_{2}F,1}(H^{a}) \,\Big)\,  t^d\,.
\end{align*}
As in (\ref{5.2})--(\ref{F-trr}), this simplifies  after noting that
(i) $(S+dF)^{2}=-n+2d$, (ii) the canonical class $K$ of $E(n)$ is
$(n-2)\,F$, and (iii) $d\,GW_{dF,1}=(2-n)\,\sigma(d)$. Thus we obtain

\begin{equation}\label{F-TRR}
     H(t)\,=\,\frac{1}{12}\,t\,F^{\prime}(t)\,  -\,
    \frac{1}{12}\,F(t)\, +\,  (2-n)\,G(t)\,F(t).
\end{equation}
On the other hand, by the same arguments for (\ref{F-sum}) and
the sum formula (\ref{FVS}) yields
\begin{equation}\label{F-SUM}
    H(t)\,=\,-\frac{1}{12}\,F(t)\  +\  2\,G(t)\,F(t).
\end{equation}
Equating (\ref{F-TRR}) and (\ref{F-SUM})   gives
the ODE
\begin{equation}\label{S9-ODE}
    \frac{1}{12}\,t\,F^{\prime}(t)\,  =\,
    n\,G(t)\,F(t)
\end{equation}
with the initial condition $F(0)=GW^{{{\mathcal H}}}_{S,0}=1$. The 
solution is the beautiful formula
\begin{equation}\label{c-E(n)}
    F(t)\,=\,\prod_{d\geq 1}
    \left(\frac{1}{1-t^{d}}\right)^{12n}
\end{equation}
which generalizes the (\ref{E:Main}) for $E(1)$. 

 Another 
application of the sum formula gives the corresponding generating 
function for higher genus curves  in $E(n)$;  see \cite{l2} for 
details.  Bryan and Leung obtain exactly  the same formulas by
using 
a different set of family GW invariants of $E(n)$ (this extends the 
$E(1)$ and $E(2)$ cases done  in \cite{bl1}).

\medskip

The exponent $12n$ in (\ref{c-E(n)}) is the Euler characteristic of 
$E(n)$.  Consequently, (\ref{c-E(n)}) is exactly the conjectured 
generating function of the rational curves in an $E(n)$ surface 
(see\cite{go}).  However, we do not know whether this count is {\em 
enumerative}, that is, whether it it the same as an actual count of 
holomorphic curves in a generic K\"{a}hler $E(n)$.

This issue arises because we have moved to a generic almost complex 
structure $J$, where we can apply symplectic cut-and-paste methods to 
find counts of curves.  We can return to the K\"{a}hler  case by 
taking a limit of generic $J$ that approach a K\"{a}hler $J_0$. 
Under that limit holomorphic curves converge to holomorphic curves, 
but degeneracies may appear:  the limit curves might be multiply 
covered, and $J_0$ may admit families of curves when only a discrete 
collection of curves is expected.  When such degeneracies are present 
it is not all clear when is meant by an enumeration of the curves. 
But this is an issue of algebraic geometry, quite separate from the 
symplectic invariant represented by the count of equation 
(\ref{c-E(n)}).

\vskip 1cm

%%%%%%%%%%%%%%%%%%%%%%%%%%%%%%%%%%%%%%%%%%%%%%%%%%%%%%%%%%%%%%%%%%%%%%%%%%%%%
%%%%%%%%%%%%%%%% Section 10  %%%%%%%%%%%%%%%%%%%%%%%%%%%%%%%%%%%%%%%%%%%%%
%%%%%%%%%%%%%%%%%%%%%%%%%%%%%%%%%%%%%%%%%%%%%%%%%%%%%%%%%%%%%%%%%%%%%%%%%%%%%

\setcounter{equation}{0}
\section{The Yau-Zaslow  Conjecture}
\label{section10}
\bigskip

About a decade ago Yau and Zaslow \cite{yz} conjectured a formula
counting rational curves in K3 surfaces.   Let $N_{d}$ be the number
of rational curves (i.e. curves whose domain has geometric genus $0$) 
in K3 surfaces
that represent some homology class $A$ with $A^{2}=2d-2$. The 
conjecture is that the
generating function for $N_{d}$ is
\begin{equation}\label{YZ}
    \sum_{d\geq 0} N_{d}\,t^{d}\,=\,
    \prod_{d\geq 1}\left(\frac{1}{1-t^{d}}\right)^{24}.
\end{equation}
This is exactly the formula obtained in the previous section for 
curves in $E(2)=K3$ in the classes $S+dF$.  The Yau-Zaslow 
Conjecture asserts that the same formula counts {\em all} rational 
curves in a generic algebraic K3 surface.

Bryan and Leung \cite{bl1} used their ``twistor family invariants'' 
to prove (\ref{YZ}) for the cases where $A$ is a primitive class.  An 
independent proof, also for primitive classes,  was given by one of 
us (J.  Lee) using symplectic  methods (\cite{l2}).  In fact,  we 
have just described that proof:
  equations
(\ref{BL}) and (\ref{c-E(n)}) directly imply  (\ref{YZ}) for 
primitive classes in $E(2)=K3$.  The key point is that for primitive 
classes the count of holomorphic {\em maps}, given by the 
family GW 
invariants, is the same as the count of holomorphic {\em curves}, 
which are counted by the numbers $N_d$.

Verifying the Yau-Zaslow  Conjecture for multiple classes has turned 
out to be much harder.  Any approach must deal with the fact that the 
counts of curves and maps are no longer the same because maps can 
multiply cover their images. Using
algebraic geometry, A. Gathmann \cite{ga} worked out some special
cases by hand, and Jun Li \cite{li} outlined an approach.  Very 
recently, the first author and N.C. Leung have extended the approach 
of the previous section to prove the Yau-Zaslow formula  (\ref{YZ}) 
for classes that
are {\em twice} a primitive class.  We will explain that result in 
this section.  The key theorem
is the following.

\begin{theorem}[\cite{ll}]
%\label{T:Main}
Suppose that $X$ is a K3 surface and $A\in H_{2}\big(X;{\Bbb Z}\big)$ is twice
a primitive class. Then the genus $g=0$ family GW-invariant of
$X$ for the class $A$ is given by
\begin{equation}
    GW_{A,0}^{{\mathcal H}}\,=\,  GW_{B,0}^{{\mathcal H}}\,+\,
    \Big(\frac{1}{2}\Big)^{3}\,GW_{A/2,0}^{{\mathcal H}}  \label{10.2}
\end{equation}
where $B$ is any primitive class with $B^{2}=A^{2}$. \end{theorem}

\medskip

 This theorem implies the  Yau-Zaslow formula (\ref{YZ}) 
for the classes that are twice a primitive class  by the following 
reasoning.   First, since we already know the  Yau-Zaslow conjecture 
for primitive classes, 
 the invariant $GW_{B,0}^{{\mathcal H}}$ 
equals to $N_{d'}$ where $d'$ is defined by $B^{2}=A^{2}=2d'-2$. 
Similarly, since $A/2$ is a primitive class, the number of curves 
representing $A/2$ is $GW^{{\mathcal H}}_{A/2,0}=N_d$ where $4(2d-2)=2d'-2$.  As 
explained by Gathmann \cite{ga}, the double covers of those 
$A/2$-curves contribute
$$
\Big(\frac{1}{2}\Big)^{3}\,GW_{A/2,0}^{{\mathcal H}} 
$$
 to the 
invariant $GW^{{\mathcal H}}_{A,0}$. Thus (\ref{10.2}) is equivalent 
to the fact that the count of (reduced) rational curves representing 
the classes $A$ and $B$ are the same --- both are the number given by formula 
(\ref{YZ}).

\medskip

\begin{proof}  The approach is the same as above:
relate the $g=1$ TRR formula and the symplectic sum formula for 
$E(2)=K3$. This is done in detail in \cite{ll}.   Here is an outline.

First note that by (\ref{BL}) we can assume that $A/2$ has the form 
$S+dF$ and $B=S+(4d-3)F$; then $A^2=B^2=8(d-1)$.  Thus  it
suffices to show that
\begin{equation*}
    GW^{{{\mathcal H}}}_{2(S+dF),0}\,= \,GW^{{{\mathcal H}}}_{S+(4d-3)F,0}\,+\,
    \Big(\frac{1}{2}\Big)^{3}\,GW_{S+dF,0}^{{{\mathcal H}}}.
\end{equation*}
This follows from

\begin{equation}\label{main-1}
    \sum GW^{{{\mathcal H}}}_{2S+dF,0}\,t^{d}\,-\,
    \sum GW^{{{\mathcal H}}}_{S+(2d-3)F,0}\,t^{d}\,=\,
    \Big(\frac{1}{2}\Big)^{3}\,
    \sum GW_{S+dF,0}^{{{\mathcal H}}}\,t^{2d}.
\end{equation}

\smallskip

\noindent Note that  if $d$ is odd  both classes $2S+dF$ and $S+(2d-3)F$ are 
primitive with the same square and hence by (\ref{BL}) all odd terms 
in the left hand side  of (\ref{main-1}) vanish.  Thus is suffices to show (\ref{main-1}).  That requires several steps.

\bigskip

We first derive a relation among the family invariants for the classes $2S+dF$,
$d\geq 0$. Fix a smooth fiber $V$ of $E(2)\to S^{2}$ and let 
    $$GW^{V}_{2S+dF,1,(2)}(C_{F})$$
denote  the genus 1 relative family invariants of $(E(2), V)$ 
 and  with multiplicity vector (2).  This is a count of stable $J_\alpha$-holomorphic maps $f$ from an elliptic curve with one marked point $x$ so that $f$ contacts $V$ at $f(x)$ with multiplicity two and whose image has no components in $V$. Introduce 
generating functions
\begin{equation*}
    M_{g}(\,\cdot\,)\,=\,\sum\,GW^{{{\mathcal H}}}_{2S+dF,g}(\,\cdot\,)\,t^{d}
  \qquad\mbox{and}\qquad
    M^{V}_{1,(2)}\,=\,\sum\,GW^{V}_{2S+dF,1,(2)}(C_{F})\,t^{d}.
\end{equation*}
The genus 1 TRR formula combined with splitting formulas gives
\begin{equation}\label{ET}
    3M_{1}\big(\tau(F)\big)\,=\,
t\,M_{0}^{\prime}\,-\,2M_{0}.
\end{equation}
while the sum formula for the symplectic sum $E(2)=E(2)\#_{V}E(0)$
yields
\begin{equation}
    M_{1}\big(\tau(F)\big)\,=\, M^{V}_{1,(2)} \,+\, 4\,G_{2}\,M_{0},
    \label{ES1}
    \end{equation}
    and
    \begin{equation}
    M_{2}\big(\tau(F),\tau(F)\big)\,-2M_{1}\big(pt\big)\, =\,
    20G_{2}\,M^{V}_{1,(2)}\, +\,\left(16G^{2}_{2}\,
    +\,8t\,G^{\prime}_{2}\right)M_{0}.
     \label{ES2}
    \end{equation}
where $G_{2}(t)$ is the Eisenstein series with weight 2, namely
   $$G_{2}(t)\,=\,-\frac{1}{24}\,+\,G(t)\,=\,
     -\frac{1}{24}\,+\,\sum_{d\geq 0}\sigma(d)\,t^{d}.$$
(The coefficients of (\ref{ES1}) and (\ref{ES2}) consist of
relative invariants of $E(0)$ that can be computed by several 
applications of genus 0 and 1
TRR formulas and the sum formula for the fiber sum
$E(0)=E(0)\#E(0)$).  Together,  equations (\ref{ET}),
(\ref{ES1}), and (\ref{ES2}) give
\setlength\multlinegap{.7cm}
\begin{multline}
\label{ODE-1}
    3M_{2}\big(\tau(F),\tau(F)\big)\,-6M_{1}\big(pt\big)\\
    \,  =\, 20t\,G_{2}M^{\prime}_{0}\,-\, \left(192G_{2}^{2}\,+\, 40
    G_{2}\,-\,24t\,G_{2}^{\prime}\right)M_{0}\,.
\end{multline}

\medskip

We can repeat this argument to get a relation  of invariants for the primitive
classes $S+(2d-3)F$.  Note that we can always find an embedded symplectic submanifold
$U$ of $E(2)$ that represents  $2F$:  take a map from an elliptic curve that double covers a fiber and perturb slightly; the perturbed map will remain symplectic, and will be an embedding by the adjunction formula since $F\cdot F=0$.  Fix such a $U$ and write
    $$GW^{U}_{S+(2d-3)F,1,(2)}(C_{2F})$$
for the genus 1 relative invariants of $E(2)$ relative to $U$  and
multiplicity vector (2). We can then build two more generating functions:
\begin{equation*}
    P_{g}(\,\cdot\,)\,=\,\sum\,GW^{{{\mathcal 
H}}}_{S+(2d-3)F,0}(\,\cdot\,)\,t^{d}
    \,,\ \ \
    P^{U}_{1,(2)}\,=\,\sum\,GW^{U}_{S+(2d-3)F,1}(C_{2F})\,t^{d}.
\end{equation*}
For these,  the genus $1$ TRR formula gives a formula like (\ref{ET})  and
  the sum formula for the symplectic sum $E(2)=E(2)\#_{U}E(0)$
gives formulas like (\ref{ES1}) and (\ref{ES2}); together those give
\setlength\multlinegap{.7cm}
\begin{multline}
  \label{ODE-2}
    3P_{2}\big(\tau(F),\tau(F)\big)\,-6\,P_{1}\big(pt\big)\\
    \,  =\, 20t\, G_{2}P^{\prime}_{0}\,-\, \left(192\,G_{2}^{2}\,+\, 40
    \,G_{2}\,-\,24\,t\,G_{2}^{\prime}\right)P_{0}\,.
\end{multline}
Note that (\ref{ODE-1}) and (\ref{ODE-2}) have the same coefficients; 
that is true because all
coefficients in  the TRR and the sum formulas
depend only on the topological quantities
\begin{equation*}
    (2S+dF)^{2}\,=\,(S+(2d-3)F)^{2},\ \ \
    (2S+dF)\cdot F\,=\,(S+(2d-3)F)\cdot 2F,\ \ \
    F^{2}\,=\,(2F)^{2}.
\end{equation*}
Hence (\ref{ODE-1}) and (\ref{ODE-2}) give
\setlength\multlinegap{.4cm}
\begin{multline}
3\big[M_{2}\big(\tau(F),\tau(F)\big) -
           P_{2}\big(\tau(2F),\tau(2F)\big)\big]
    - 6\big[M_{1}\big(pt\big) -  P_{1}\big(pt\big)\big] \\
    \, =\, 20 G_{2}t\,\big[M_{0} - P_{0}\,\big]^{\prime} -
    \left(192G_{2}^{2} + 40G_{2} -  24t\,G_{2}^{\prime}\right)
    \big[M_{0} - P_{0}\big]
     \label{ODE-3}
\end{multline}

\medskip

Next observe that because $2S+dF$ and $S+(2d-3)F$ are primitive 
classes  with the
same square when $d$ is odd, it follows from (\ref{BL}) that the
generating functions  $M_{1}\big(pt\big) - P_{1}\big(pt\big)$
and $M_{0}-P_{0}$ have no odd terms. One can also show that
   $M_{2}\big(\tau(F),\tau(F)\big) -
           P_{2}\big(\tau(2F),\tau(2F)\big)$
has no odd terms (see \cite{ll2}). Therefore, comparing odd terms
in both sides of (\ref{ODE-3}) gives the first order ODE
\begin{equation}  \label{ODE-4}
    0 =\, 20G_{o}t\,
    \big(M_{0}-P_{0}\big)^{\prime} -
    \left(384\,G_{e}\,G_{o}\,+\,40\,G_{o} - 
24t\,G_{o}^{\prime}\right)\big(M_{0}-P_{0}\big)\,
\end{equation}
where $G_{e}(t)$ (resp. $G_{o}(t)$) is the sum of all even (resp. odd)
terms of $G_{2}(t)$.  The initial condition
of this ODE is given by the well-known multiple map contribution
\cite{ga}
    $$(M_{0}-P_{0})(0)\,=\,GW^{{{\mathcal H}}}_{2S,0}\,=\,
      \left(\frac{1}{2}\right)^{3}. $$

\smallskip
On the other hand, it follows from the differential equation
(\ref{S9-ODE}) that
\begin{equation*}
    t\,\frac{d}{dt}F(t^{2})\,=\,
    2\,t^{2}\,F^{\prime}(t^{2})\,=\,48\,G(t^{2})\,F(t^{2})\,=\,
    48\,G_{2}(t^{2})\,F(t^{2})\,+\,2\,F(t^{2}).
\end{equation*}
Using this, one can check that $F(t^{2})$ satisfies the same ODE
(\ref{ODE-4}) if
\begin{equation}  \label{rel-qmod}
-4\,t^{2}\,G_{2}^{\prime}(t^{2})\,+\,32\,G_{2}^{2}(t^{2})\,-\,
40\,G_{2}(t)\,G_{2}(t^{2})\,+\,8\,G_{2}^{2}(t)\,-\,t\,G_{2}^{\prime}(t)
\,=\,0.
\end{equation}
Under the substitution $t=e^{2\pi i z}$ the left hand side of (\ref{rel-qmod}) is a modular form of level
2 with weight 4 and the space of those modular forms is a 2
dimensional vector space with well-known bases \cite{k}. Thus, one
can prove (\ref{rel-qmod}) by showing the first two terms in the expansion of the
left hand side of (\ref{rel-qmod}) vanish. Using the initial
condition $F(0)=GW^{{{\mathcal H}}}_{S,0}=1$, we can thus conclude that
\begin{equation*}
\label{mul-con}
M_{0}(t)-P_{0}(t)\,=\,\Big(\frac{1}{2}\Big)^{3}\,F(t^{2}).
\end{equation*}
By the definition of $M_{0}(t)$ and $P_{0}(t)$, this is exactly the 
desired formula
(\ref{main-1}). \end{proof}

\bigskip

It would be interesting to extend this approach to higher multiple 
classes.  Unfortunately that is not straightforward.  In the above 
calculation,    the right-hand side of the sum formula  (\ref{ES2}) 
fortuitously contains only  genus 0 and 1 invariants of $E(2)$.  That 
is a consequence of the vanishing of the 
genus 0 relative invariants 
of $E(0)$  constrained to contact three copies of $F$ to order 2.
That vanishing does not occur for the classes $aS+dF$ with $a\geq 3$. 
That leads to additional terms in   (\ref{ES2})  for  higher 
multiple classes, and dealing with those terms seems to require new 
ideas.

    \vskip1cm

{\small

%%%%%%%%%%%%%%%%%%%%%%%%%%%%%%%%%%%%%%%%%%%%%%%%%%%%%%%%%%%%%%%%%%%%%%%%%%%%%%%

}


\begin{thebibliography}{99}

\bibitem[BL1]{bl1} J. Bryan and N. C. Leung,
                      \textit{The enumerative geometry of K3 surfaces
                              and modular forms}, J. Amer. Math. Soc.
                         \textbf{13} (2000), 371--410.

\bibitem[BL2]{bl2} J. Bryan and N. C. Leung, \textit{Genenerating functions
for the number of curves on abelian surfaces}, Duke Math. J.
\textbf{99} (1999), no. 2, 311--328.

\bibitem[BL3]{bl3} J. Bryan and N. C. Leung, \textit{Counting curves on
irrational surfaces}, Surveys in differential geometry:
differential geometry inspired by string theory, 313--339, Surv.
Diff. Geom., 5, Int. Press, Boston, MA, 1999.



\bibitem[BF]{bf} K. Behrend and B. Fantechi, In Preparation.



\bibitem[Ga]{ga}  A. Gathmann, \textit{The number of plane conics 5-fold
tangent to a given curve}, preprint, math.AG/0202002.



\bibitem[Ge]{ge}  E. Getzler, \textit{Topological recursion relations
in genus 2}
, In ''Integrable systems and algebraic geometry (Kobe/Kyoto,
1997).'' World Sci. Publishing, River Edge, NJ, 198, pp 73--106.


\bibitem[G\"{o}]{go} L. G\"{o}ttsche, \textit{A conjectural 
generating function for numbers of curves on surfaces},  Comm. Math. 
Phys. \textbf{196} (1998), no. 3, 523--533.




\bibitem[IP1]{ip1} E. Ionel and T.  Parker,  \textit{Gromov-Witten Invariants
of Symplectic Sums}, Math. Res. Lett.,  \textbf{5}(1998), 563--576.

\bibitem[IP2]{ip2} E. Ionel and T.  Parker,  \textit{Relative
Gromov-Witten Invariants}, Annals of Math. \textbf{157} (2003), 45--96.


\bibitem[IP3]{ip3} E. Ionel and T.  Parker,  \textit{The Symplectic Sum
Formula for Gromov-Witten Invariants}, Annals of  Math., \textbf{159} (2004), 935-1025.




\bibitem[K]{k}  A.W. Knapp, \textit{Elliptic curves}, Princeton University
Press, 1992.


\bibitem[L1]{l1}  J. Lee, \textit{Family Gromov-Witten Invariants for
K\"{a}hler Surfaces},  Duke Math. J. \textbf{ 123} (2004), no 1, 209--233.

\bibitem[L2]{l2}  J. Lee, \textit{Counting Curves in Elliptic Surfaces by
Symplectic Methods}, preprint, math. SG/0307358.


\bibitem[LL1]{ll}  J. Lee and N.C. Leung, \textit{Yau-Zaslow formula on K3
surfaces for non-primitive classes}, preprint, math. SG/0404537.


\bibitem[LL2]{ll2}  J. Lee and N.C. Leung, \textit{Counting elliptic curves in K3
surfaces}, preprint, math. SG/0405041.




\bibitem[Li]{li}  J. Li, \textit{A note on enumerating rational curves in a K3
surface}, Geometry and nonlinear partial differential equations
(Hangzhou, 2001), 53--62, AMS/IP Stud. Adv. Math. \textbf{29}, Amer.
Math. Soc., Providence, RI, 2002.


\bibitem[LT]{LT} J. Li and G. Tian,  \textit{Virtual moduli cycles and
Gromov-Witten invariants of general symplectic manifolds}, Topics
in symplectic $4$-manifolds (Irvine, CA, 1996), 47--83, First Int.
Press Lect. Ser., I, International Press, Cambridge, MA, 1998.


\bibitem[P]{p} T. Parker, \textit{Compactified moduli spaces of
pseudo-holomorphic curves} Mirror symmetry, III (Montreal, PQ,
1995), 77--113, AMS/IP Stud. Adv. Math., 10, Amer. Math. Soc.,
Providence, RI, 1999.

\bibitem[PW]{pw} T. Parker and J. Wolfson,
\textit{Pseudo-holomorphic maps and bubble trees}, Jour. Geometric
Analysis, \textbf{3} (1993) 63--98.


\bibitem[RT1]{rt1} Y. Ruan and G. Tian, \textit{A mathematical theory
of quantum
cohomology}, J. Differential Geom. \textbf{42} (1995), 259--367.

\bibitem[RT2]{rt2}  Y. Ruan and G. Tian, \textit{Higher genus symplectic
invariants and sigma models coupled with gravity},  Invent. Math.
\textbf{130} (1997), 455--516.

\bibitem[IS]{is} S. Ivashkovich and V. Shevchishin,
\textit{Gromov compactness theorem for $J$-complex curves with
boundary}, Internat. Math. Res. Notices \textbf{22} (2000),
1167--1206.

\bibitem[YZ]{yz} S.T. Yau and E. Zaslow, \textit{BPS States, String Duality,
and Nodal Curves on K3}, Nuclear Phys. B \textbf{471} (1996), 503--512.


\end{thebibliography}
\end{document}